\documentclass[12pt]{amsart}
\usepackage{pgfplots}
\usepackage{tikz}

\usetikzlibrary{automata}
\usetikzlibrary{positioning}
\usepackage{amsmath, amssymb} 
\usepackage{amsthm}
\usepackage{graphicx}
\usepackage[shortlabels]{enumitem}

\usepackage{amsfonts,geometry,hyperref}

\newtheorem{conj}{Conjecture}[section]

\newtheorem{theorem}{Theorem}[section]
\newtheorem{prop}[theorem]{Proposition}
\newtheorem{lemma}[theorem]{Lemma}
\newtheorem{cor}[theorem]{Corollary}

\DeclareMathOperator{\ESCom}{ESCom}
\DeclareMathOperator{\CSCom}{CSCom}
\DeclareMathOperator{\OSCom}{OSCom}
\DeclareMathOperator{\OrdClass}{OC}
\DeclareMathOperator{\CClass}{Cl}

\begin{document}
\title[Super commuting graphs of  finite groups and their Zagreb indices]{Super commuting graphs of  finite groups and their Zagreb indices}

\author[S. Das and R. K. Nath]{Shrabani Das and Rajat Kanti Nath}

\address{S. Das, Department of Mathematical Sciences, Tezpur University, Napaam-784028, Sonitpur, Assam, India.\newline
Department of Mathematics, Sibsagar University, Joysagar-785665, Sibsagar, Assam, India.}

\email{shrabanidas904@gmail.com}

\address{R. K. Nath, Department of Mathematical Sciences, Tezpur University, Napaam-784028, Sonitpur, Assam, India.} 
\email{ rajatkantinath@yahoo.com}

\begin{abstract}
Let $B$ be an equivalence relation defined on a finite group $G$. The $B$ super commuting graph of $G$ is a graph whose vertex set is $G$ and two distinct vertices $g$ and $h$ are adjacent if either $[g] = [h]$ or there exist $g' \in [g]$ and $h' \in [h]$ such that $g'$ commutes with $h'$, where  $[g]$ is the  $B$-equivalence class of $g \in G$.	Considering  $B$ as the equality, conjugacy and same order relations on $G$, in this article, we discuss the graph structures of equality/conjugacy/order super commuting graphs of certain well-known families of non-abelian groups viz. dihedral groups, generalized quaternion groups, semidihedral groups, quasidihedral groups, the groups $U_{6n}, V_{8n}, M_{2mn}$ etc. 
Further, we compute the Zagreb indices of these graphs and show that they satisfy the Hansen--Vuki{\v{c}}evi{\'c} conjecture.     
\end{abstract}

\thanks{ }

\maketitle
\section{Introduction}
For quite some time, graphs have been defined on finite groups and their properties have been studied through the properties of groups and vice versa. Finite groups are also characterized subject to the satisfaction of certain  graph theoretic properties of various graphs defined on the groups.  Many times these graphs produce examples/counter examples for various Conjectures/Problems in graph theory.
The commuting graph,  non-commuting graph,  generating graph,  power graph, nilpotent graph, solvable graph etc. are some well-known graphs defined on groups. Another class of graphs has been defined on groups by considering the vertex set as the set of  conjugacy classes and two conjugacy classes are adjacent  if there exist  elements in the respective conjugacy  classes satisfying some group property. Commuting/nilpotent/solvable conjugacy class graph are examples of  such graphs. By mixing these two types of graphs defined on groups Arunkumar et al. \cite{GA-PJC-RKN-LS-2022} have introduced the notion of super graphs of groups,  recently  in 2022.


 Let $G$ be a finite group and $B$  an equivalence relation defined on it. We write  $[g]$ to denote the  $B$-equivalence class of $g \in G$.
 The $B$ super commuting graph on $G$ is a graph whose vertex set is $G$ and two distinct vertices $g$ and $h$ are adjacent if either $[g] = [h]$ or there exist $g' \in [g]$ and $h' \in [h]$ such that $g'$ commutes with $h'$.
 Considering $B$ as the equality,  conjugacy and same order relations on the elements of $G$, we have equality super commuting graph (denoted by $\ESCom(G)$),  conjugacy super commuting graph (denoted by $\CSCom(G)$) and order super commuting graph (denoted by $\OSCom(G)$) respectively. These classes of graphs on finite groups were introduced in  \cite{GA-PJC-RKN-LS-2022} and studied their properties in \cite{GA-PJC-RKN-LS-2022,DMP-24}.  The Wiener index of equality/conjugacy super commuting   graphs for dihedral and generalized quaternion groups have been computed in \cite{GA-PJC-RKN-2024} while their spectrum in \cite{GA-PJC-RKN-GB-2024}. The spectrum and Laplacian spectrum of order super commuting graphs of the dihedral groups, generalized quaternion groups and non-abelian groups of order $pq$ (where $p$ and $q$ are primes) are computed in \cite{DMP-24-arxiv}.  Note that  equality/conjugacy/order super commuting graphs of those groups were realized through generalized composition of graphs while computing various spectra and Wiener index. In this paper, we realize the graph structures of equality/conjugacy/order super commuting graphs of  the dihedral groups, generalized quaternion groups, semidihedral groups, quasidihedral groups, the groups $U_{6n}, V_{8n}$ and $M_{2mn}$  in terms of join and union of certain complete graphs. These realizations reduce the complexity in computing various graph theoretic parameters. We demonstrate this by computing certain topological indices viz. first and second Zagreb indices of these graphs and show that  they satisfy the Hansen--Vuki{\v{c}}evi{\'c} conjecture.  It is worth mentioning that in \cite{DEN-25,DSN-25,DSN-25-arxiv}, Das et al. have computed first and second Zagreb indices of subgroup generating bipartite graphs, commuting/non-commuting graphs, commuting conjugacy class graphs of certain groups respectively and showed that they satisfy Hansen--Vuki{\v{c}}evi{\'c} conjecture. Therefore, it is natural to ask the question whether Hansen--Vuki{\v{c}}evi{\'c} conjecture is true for equality/conjugacy/order super commuting graphs.

  Let $\mathfrak{G}$ be the set of all graphs. A topological index  is a function $T :  \mathfrak{G} \to \mathbb{R}$ such that $T(\Gamma_1) = T(\Gamma_2)$ whenever the graphs $\Gamma_1$ and   $\Gamma_2$ are isomorphic. Since 1947 numerous topological indices have been defined by using different parameters of graphs. Among them, Zagreb index is a very popular degree-based topological index. It  was introduced by Gutman and Trinajsti{\'c} \cite{Gut-Trin-72} in 1972. Zagreb index was initially used in examining the dependence of total $\pi$-electron energy on molecular structure. It was also used in studying molecular complexity, chirality, ZE-isomerism and heterosystems etc. \cite{Z-index-30y-2003}. We can find the mathematical properties of Zagreb index in the survey \cite{Gut-Das-2004}. 

   Let $\Gamma$ be a simple undirected graph with vertex set $v(\Gamma)$ and edge set $e(\Gamma)$. The first and second Zagreb indices of $\Gamma$, denoted by $M_{1}(\Gamma)$ and $M_{2}(\Gamma)$ respectively, are defined as 
\[
M_{1}(\Gamma) = \sum\limits_{v \in v(\Gamma)} \deg(v)^{2}  \text{ and }  M_{2}(\Gamma) = \sum\limits_{uv \in e(\Gamma)} \deg(u)\deg(v),
\]
where $ \deg(v) $ is the number of edges incident on $ v $
(called degree of $v$). Comparing first and second Zagreb indices, Hansen and Vuki{\v{c}}evi{\'c} \cite{hansen2007comparing} posed the following   conjecture in  2007.
\begin{conj}\label{Conj}
	(Hansen--Vuki{\v{c}}evi{\'c} Conjecture) For any simple finite graph $\mathcal{G}$, 
	\begin{equation}\label{Conj-eq}
		\dfrac{M_{2}(\mathcal{G})}{\vert e(\mathcal{G}) \vert} \geq \dfrac{M_{1}(\mathcal{G})}{\vert v(\mathcal{G}) \vert} .
	\end{equation}
\end{conj}
It was shown in \cite{hansen2007comparing} that the conjecture is not true  if $\Gamma = K_{1, 5} \sqcup K_3$. However, Hansen and Vuki{\v{c}}evi{\'c} \cite{hansen2007comparing}  showed that  Conjecture \ref{Conj}  holds for chemical graphs. 
Further studies in \cite{vukicevic2007comparing, liu2008conjecture} showed that the conjecture holds for trees and connected unicyclic graphs. A survey on comparing Zagreb indices can be found in \cite{Liu-You-2011}.


\textbf{Notations:} We write $\vee$ to denote the join of two graphs, $\sqcup$/$\bigsqcup$  to denote the union of two disjoint sets or graphs. For any two vertices $x$ and $y$ of a graph $\Gamma$, if $x$ is adjacent to $y$ then we write  $x \sim y$ otherwise $x \nsim y$. The disjoint union of $m$ copies of $K_n$ is denoted by $mK_n$. Let $\mathcal{H}$ be any graph with vertex set $v(\mathcal{H})$.  For any subset $V'$ of $v(\mathcal{H})$ we write  $\mathcal{H}[V']$ to denote the induced subgraph of $\mathcal{H}$ induced by $V'$.

\section{Realizing super graphs}
 Schwenk  \cite{AJS-2006} introduced the following notion of generalized composition of graphs, denoted it by $\mathcal{H}[\Gamma_1, \Gamma_2, \ldots, \Gamma_k]$, where $\mathcal{H}$ is a graph with vertex set $v(\mathcal{H})=\{1, 2, \ldots, k\}$ and  $\Gamma_1, \Gamma_2, \ldots, \Gamma_k$ are  graphs with  vertex sets $v(\Gamma_i)=\{ v^{1}_i, \ldots, v^{n_i}_i\}$ for $1 \leq i \leq k$.
The  generalized composition $\mathcal{G}= \mathcal{H}[\Gamma_1, \Gamma_2, \ldots, \Gamma_k]$  is a graph whose  vertex set is $v(\Gamma_1) \sqcup \cdots \sqcup v(\Gamma_k)$ and two vertices $v^{p}_i$ and $v^{q}_j$ of $\mathcal{G}$ are adjacent if either of the conditions hold:
\begin{enumerate}
\item $i=j$ and $v^{p}_i$, $v^{q}_j$ are adjacent vertices in $\Gamma_i$.
\item $i \neq j$ and $i$, $j$ are adjacent in $\mathcal{H}$.
\end{enumerate}
Note that this notation for generalized composition of graphs  should not be confused with the notation for induced subgraph. In that one we have a subset of vertices inside the brackets, and in this one we have graphs.
 
In \cite{GA-PJC-RKN-2024}, equality and conjugacy super commuting graphs of the dihedral groups $D_{2n} = \langle a, b :  a^n= b^2 = 1, bab^{-1} = a^{-1} \rangle$ (for $n\geq 3$) and  generalized quaternion groups $Q_{4n}=\langle a, b: a^{2n}=1, a^n=b^2, bab^{-1}=a^{-1} \rangle$ (for $n\geq 2$) have been realized as generalized composition of certain graphs. In this section, we realize those graphs as join of certain graphs.
Further, we shall realize order super commuting graph of these groups. In our study, we shall also consider the quasidihedral groups $QD_{2^n} = \langle a, b : a^{2^{n-1}} = b^2 = 1, bab^{-1} = a^{2^{n-2}-1} \rangle$ (for $n\geq 3$), the semidihedral group $SD_{8n} = \langle a, b: a^{4n}=b^{2}=1, bab^{-1}=a^{2n-1} \rangle$ (for $n \geq 2$) and the  groups $V_{8n}=\langle a, b: a^{2n}=b^{4}=1, ba=a^{-1}b^{-1}, b^{-1}a=a^{-1}b\rangle$, $U_{6n} = \langle a, b: a^{2n}=b^3=1, a^{-1}ba=b^{-1} \rangle$ and $M_{2mn}= \langle a, b: a^{m}=b^{2n}=1, bab^{-1}=a^{-1} \rangle$ (for $m \geq 3$ but not equal to $4$).

\subsection{Equality super commuting graphs}
The equality super commuting graph of a group $G$  (denoted by $\ESCom(G)$) is a graph with vertex set $v(\ESCom(G))=G$ and two distinct vertices $g$ and $h$ are adjacent if either $[g]=[h]$ or there exists $g' \in [g]$ and $h' \in [h]$ such that $g'h'=h'g'$, where $[g]$ denotes the  equivalence classes with respect to the equivalence relation on $G$ namely `equality'. Since a group has distinct elements,  $[g] = \{g\}$  for all $g \in G$ and $[g] \neq [h]$ if $g \neq h$. Thus, $\ESCom(G)$ is a graph whose vertex set is $\{\{g\}: g \in G\}$ and two distinct vertices $\{g\}$ and $\{h\}$ are adjacent if $gh=hg$. It follows that $\ESCom(G) \cong C(G)$, where $C(G)$ is the commuting graph of $G$ which is a graph with vertex set $G$ and where two distinct vertices are adjacent if they commute. 
%
Since central elements of a group commute with every other element of the group, we have
\[
\ESCom(G) \cong K_{|Z(G)|} \vee C(G \setminus Z(G)),
\]
where  $Z(G)$ is the center of $G$ and $C(G \setminus Z(G))$ is the induced subgraph $C(G)[G \setminus Z(G)]$ of $C(G)$ induced by $G \setminus Z(G)$.
Note that the centers and  subgraphs of $C(G)$ induced by the non-central elements of the groups mentioned above are well known (see \cite{Dut-Nath-2017,JD-RKN-2017,mirzargar2012some}). Therefore, we have the following result.
\begin{prop}\label{ESCom(G)}
\begin{enumerate}
\item $
	\ESCom(D_{2n})\cong \begin{cases}
		K_1 \vee (nK_1 \sqcup K_{n-1}), & \text{if } $n$ \text{ is odd} \\
		K_2 \vee (\frac{n}{2}K_2 \sqcup K_{n-2}), & \text{  if } $n$ \text{ is even}.
	\end{cases}
$

\vspace{.2cm}

\item $\ESCom(Q_{4n}) \cong K_2 \vee \left(n K_2 \sqcup K_{2n-2}\right)$.

\item $\ESCom(QD_{2^n})\cong K_2 \vee (K_{2^{n-1}-2} \sqcup 2^{n-2}K_2)$.

\item $\ESCom(V_{8n})\cong \begin{cases}
	K_2 \vee (K_{4n-2} \sqcup 2nK_2), & \text{  if } $n$ \text{ is odd} \\
	K_4 \vee (K_{4n-4} \sqcup nK_4), & \text{  if } $n$ \text{ is even}.
\end{cases}$

\item $\ESCom(SD_{8n}) \cong \begin{cases}
 K_4 \vee (K_{4n-4} \sqcup nK_4), & \text{  if } n \text{ is odd} \\
 K_2 \vee (K_{4n-2} \sqcup 2nK_2), & \text{  if } n \text{ is even}.
 \end{cases}$

\item $\ESCom(U_{6n}) \cong K_n \vee (K_{2n} \sqcup 3K_n)$.

\item $\ESCom(M_{2mn})\cong \begin{cases}
	K_n \vee (K_{(m-1)n} \sqcup mK_n), & \text{  if } $m$ \text{ is odd} \\
	K_{2n} \vee (K_{(\frac{m}{2}-1)2n} \sqcup \frac{m}{2}K_{2n}), & \text{  if } $m$ \text{ is even}.
\end{cases}$
\end{enumerate}
\end{prop}

 \subsection{Conjugacy super commuting graphs}\label{CScom}
 The conjugacy super commuting graph of a group $G$, denoted by $\CSCom(G)$, is a graph with vertex set $v(\CSCom(G))=G$ and two distinct vertices $g$ and $h$ are adjacent if either $[g]=[h]$ or there exists $g' \in [g]$ and $h' \in [h]$ such that $g'h'=h'g'$, where $[g] = \CClass(g)$ is the conjugacy  class of $g$ in $G$.
In \cite{GA-PJC-RKN-2024}, Arunkumar et al. obtained the following graph structures of $\CSCom(G)$ (as generalized composition of graphs) when $G$ is $D_{2n}$ and $Q_{4n}$: 
\begin{align*}
&\CSCom(D_{2n})\\
 &\quad \quad\cong \begin{cases}
    (K_1 \vee (K_1 \sqcup K_{\frac{n-1}{2}})) [K_1, \underbrace{K_2, \ldots, K_2}_{(\frac{n-1}{2}) \text{-times}}, K_n], \text{ if } n \text{ is odd} \\
    (K_2 \vee (K_1 \sqcup K_1 \sqcup K_{\frac{n}{2}-1})) [K_1, K_1, \underbrace{K_2, \ldots, K_2}_{(\frac{n}{2}-1) \text{-times}}, K_{\frac{n}{2}}, K_{\frac{n}{2}}], \text{ if } n \text{ and } \frac{n}{2} \text{ are even} \\
    (K_2 \vee (K_2 \sqcup K_{\frac{n}{2}-1})) [K_1, K_1, \underbrace{K_2, \ldots, K_2}_{(\frac{n}{2}-1) \text{-times}}, K_{\frac{n}{2}}, K_{\frac{n}{2}}], \text{ if } n \text{ is even and } \frac{n}{2} \text{ is odd}
\end{cases}
\end{align*}
and
\[
\CSCom(Q_{4n}) \cong \begin{cases}
    (K_2 \vee (K_1 \sqcup K_1 \sqcup K_{n-1})) [K_1, K_1, \underbrace{K_2, \ldots, K_2}_{(n-1)\text{-times}}, K_n, K_n], & \text{ if } n \text{ is even} \\
    (K_2 \vee (K_2 \sqcup K_{n-1})) [K_1, K_1, \underbrace{K_2, \ldots, K_2}_{(n-1)\text{-times}}, K_n, K_n], & \text{ if } n \text{ is odd.}
\end{cases}
\]
However, we have obtained simpler structures for $\CSCom(D_{2n})$ and $\CSCom(Q_{4n})$ using the following lemma. Our method also enables to describe $\CSCom(G)$ for the other groups mentioned in Section 2.  
\begin{lemma}\label{Lemma1}
    Let $G$ be a finite non-abelian group with center $Z(G)$. If $\Gamma(G)$ is the $B$ super commuting graph of $G$ then the following hold:
    \begin{enumerate}
        \item $\Gamma(G) \cong K_{|Z(G)|} \vee \Gamma(G)[G \setminus Z(G)]$.
        \item $\Gamma(G)[[g]] = K_{|[g]|}$ for all $g \in G$.
        \item Let $A$ be the set formed by a unique representative of each $B$-equivalence class. Let $A_1, A_2 \subseteq A$ be such that $A_1 \cap A_2=\phi$. Then
        \begin{enumerate}[label=(\roman*)]
            \item If for all $x \in A_1$ and $y \in A_2$ there exist $x' \in [x]$ and $y' \in [y]$ such that $x'y'=y'x'$, we have
            \begin{align*}
                \Gamma(G)\bigg[\bigsqcup_{g \in A_1 \sqcup A_2} [g] \bigg] \cong \Gamma(G)\bigg[\bigsqcup_{g \in A_1 }[g] \bigg] \vee   
                \Gamma(G)\bigg[\bigsqcup_{g \in A_2}[g]\bigg].
                \end{align*}
            \item If for all $x \in A_1$ and $y \in A_2$ and for all $x' \in [x]$ and $y' \in [y]$ we have $x'y' \neq y'x'$, then
                \begin{align*}
                \Gamma(G)\bigg[\bigsqcup_{g \in A_1 \sqcup A_2}[g]\bigg] \cong \Gamma(G)\bigg[\bigsqcup_{g \in A_1}[g]\bigg] \sqcup   \Gamma(G)\bigg[\bigsqcup_{g \in A_2}[g] \bigg]. 
                \end{align*}
        \end{enumerate}
    \end{enumerate}
\end{lemma}
\begin{proof}
    The proof of (a) follows from the fact that central elements commutes with every element of the group. The proof of (b) and (c) follows from the definition of $B$ super commuting graph of a group. 
\end{proof}
\begin{theorem}\label{CSCom D_2n}
Let $D_{2n}  = \langle a, b :  a^n= b^2 = 1, bab^{-1} = a^{-1} \rangle$ $(n\geq 3)$ be the dihedral group. Then
\begin{align*}
        \CSCom(D_{2n}) \cong \begin{cases}
            K_1 \vee \left(K_{n-1} \sqcup K_n \right), &\text{ if } n \text{ is odd} \\
            K_2 \vee \left(K_{n-2} \sqcup 2K_{\frac{n}{2}}\right), & \text{ if } n \text{ and } \frac{n}{2} \text{ are even} \\
            K_2 \vee \left(K_{n-2} \sqcup K_n \right), & \text{ if } n \text{ is even and } \frac{n}{2} \text{ is odd.}
            \end{cases}
    \end{align*}
\end{theorem}
\begin{proof}
If $n$ is odd then  $Z(D_{2n})=\{1\}$ and $\CClass(b) =\{b, ab, \ldots, a^{n-1}b\}$ and $\CClass(a^j) =\{a^j, a^{n-j}\}$,  where $1 \leq j \leq \frac{n-1}{2}$, are the  distinct conjugacy classes of $D_{2n}$. Let $A=\{b, a, a^2, a^3, \ldots, a^{\frac{n-1}{2}}\}, A_1=\{b\}, A_2=\{a, a^2, \ldots, a^{\frac{n-1}{2}}\}$. As such $\CSCom(D_{2n})[\CClass(b)]\cong K_{n}$. Since all powers of $a$ commute with each other, therefore, $\CSCom(D_{2n})[\bigsqcup\limits_{g \in A_2} \CClass(g)] \cong \underbrace{K_2 \vee \cdots \vee K_2}_{\frac{n-1}{2}\text{-times}}=K_{n-1}$. Also, for all $x \in A_1$ and $y \in A_2$ and for all $x' \in \CClass(x)$ and $y' \in \CClass(y)$ we have $x'y' \neq y'x'$. As such, $\CSCom(D_{2n})\Big[\bigsqcup\limits_{g \in A_1 \sqcup A_2} \CClass(g)\Big] \cong K_n \sqcup K_{n-1}$. Therefore, by Lemma \ref{Lemma1}, we have 
\begin{align*}
    \CSCom(D_{2n}) \cong K_1 \vee \CSCom(D_{2n})\Big[\bigsqcup\limits_{g \in A} \CClass(g)\Big] &  \cong K_1 \vee \Big[\bigsqcup\limits_{g \in A_1 \sqcup A_2} \CClass(g)\Big] \\
    & \cong K_1 \vee (K_{n-1} \sqcup K_n).
\end{align*}
If $n$ is even then $Z(D_{2n})=\{1, a^{\frac{n}{2}})$ and $\CClass(b) = \{b, a^2b, a^4b, \ldots, a^{n-2}b\}$, $\CClass(ab)=\{ab, a^3b, \ldots, a^{n-1}b\}$ and $\CClass(a^j) = \{a^j, a^{n-j}\}$,
where $1 \leq j \leq \frac{n}{2}-1$, are the distinct conjugacy classes of $D_{2n}$. Let $A=\{b, ab, a, a^2, a^3, \ldots, a^{\frac{n}{2}-1}\}, A_1=\{b\}, A_2=\{ab\}, A_3=\{a, a^2, \ldots, a^{\frac{n}{2}-1}\}$. As such
$\CSCom(D_{2n})[\CClass(b)] \cong  K_{\frac{n}{2}} 
$, $\CSCom(D_{2n})[\CClass(ab)] \cong K_{\frac{n}{2}}$ and $\CSCom(D_{2n})[\bigsqcup\limits_{g \in A_3} \CClass(g)] \cong \underbrace{K_2 \vee \cdots \vee K_2}_{(\frac{n}{2}-1)\text{-times}}=K_{n-2}$. For all $x \in A_1 \sqcup A_2$ and $y \in A_3$ and for all $x' \in \CClass(x)$ and $y' \in \CClass(y)$ we have $x'y' \neq y'x'$. Also, if $\frac{n}{2}$ is odd, then for $b \in A_1$ and $ab \in A_2$ there exist $b \in \CClass(b)$ and $a^{\frac{n}{2}}b \in \CClass(ab)$ such that $ba^{\frac{n}{2}}b=a^{\frac{n}{2}}bb$ but if $\frac{n}{2}$ is even, then for all $x \in A_1$ and $y \in A_2$ and for all $x' \in \CClass(x)$ and $y' \in \CClass(y)$ we have $x'y' \neq y'x'$. Therefore, 
\small{\begin{align*}
    \CSCom(D_{2n})\Big[\bigsqcup_{g \in A_1 \sqcup A_2}\CClass(g)\Big] &\cong \begin{cases}
        \CSCom(D_{2n})[\CClass(b)] \vee \CSCom(D_{2n})[\CClass(ab)], \text{ if } \frac{n}{2} \text{ is odd} \\
        \CSCom(D_{2n})[\CClass(b)] \sqcup \CSCom(D_{2n})[\CClass(ab)], \text{ if } \frac{n}{2} \text{ is even}
    \end{cases} \\
    & \cong \begin{cases}
        K_{\frac{n}{2}} \vee K_{\frac{n}{2}} = K_n,  \text{ if } \frac{n}{2} \text{ is odd }\\
K_{\frac{n}{2}} \sqcup K_{\frac{n}{2}}=2K_{\frac{n}{2}},  \text{ if } \frac{n}{2} \text{ is even.}
    \end{cases}
\end{align*}}
Thus by Lemma \ref{Lemma1} and after the above discussion, we have
\begin{align*}
    \CSCom(D_{2n}) & \cong K_2 \vee \CSCom(D_{2n})\Big[ \bigsqcup_{g \in A}\CClass(g) \Big] \\
    & \cong K_2 \vee \CSCom(D_{2n})\Big[ \bigsqcup_{g \in A_1 \sqcup A_2 \sqcup A_3}\CClass(g) \Big] \\
    & \cong \begin{cases}
        K_2 \vee \left(K_{n-2} \sqcup 2K_{\frac{n}{2}}\right), & \text{ if } n \text{ and } \frac{n}{2} \text{ are even} \\
            K_2 \vee \left(K_{n-2} \sqcup K_n \right), & \text{ if } n \text{ is even and } \frac{n}{2} \text{ is odd.}
    \end{cases}
\end{align*}
\end{proof}
\begin{theorem}\label{CSCom Q_4n}
    Let $Q_{4n}=\langle a, b: a^{2n}=1, a^n=b^2, bab^{-1}=a^{-1} \rangle$ $(n \geq 2)$ be the generalized quaternion group. Then
    \begin{align*}
        \CSCom(Q_{4n}) \cong \begin{cases}
            K_2 \vee \left(K_{2n-2} \sqcup 2K_n\right), & \text{ if } n \text{ is even} \\
             K_2 \vee \left( K_{2n-2} \sqcup K_{2n}\right), & \text{ if } n \text{ is odd.} 
             \end{cases}    
    \end{align*}
\end{theorem}
\begin{proof}
    We know $Z(Q_{4n})=\{1, a^n\}$, and $\CClass(b)=\{b, a^2b, a^4b, \ldots, a^{2n-2}b\}$, $\CClass(ab)=\{ab, a^3b$, $ \ldots, a^{2n-1}b\}$ and $\CClass(a^j)=\{a^j, a^{2n-j}\}$, where $1 \leq j \leq n-1$, are the distinct conjugacy classes of $Q_{4n}$. Let $A=\{b, ab, a, a^2, a^3, \ldots, a^{n-1}\}, A_1=\{b\}, A_2=\{ab\}, A_3=\{a, a^2, \ldots, a^{n-1}\}$. As such, we have $\CSCom(Q_{4n})[\CClass(b)]$ $\cong K_n$ and $\CSCom(Q_{4n})[\CClass(ab)]$ $\cong K_n$. Since all powers of $a$ commute with each other, therefore, $\CSCom(Q_{4n})[\bigsqcup\limits_{g \in A_3} \CClass(g)]$ $\cong \underbrace{K_2 \vee \cdots \vee K_2}_{(n-1)\text{-times}}=K_{2n-2}$. For all $x \in A_1 \sqcup A_2$ and $y \in A_3$ and for all $x' \in \CClass(x)$ and $y' \in \CClass(y)$ we have $x'y' \neq y'x'$. Also, if $n$ is odd, then for $b \in A_1$ and $ab \in A_2$ there exist $b \in \CClass(b)$ and $a^nb \in \CClass(ab)$ such that $ba^nb=a^nbb$ but if $n$ is even, then for all $x \in A_1$ and $y \in A_2$ and for all $x' \in \CClass(x)$ and $y' \in \CClass(y)$ we have $x'y' \neq y'x'$. Therefore
    \small{\begin{align*}
    \CSCom(Q_{4n})\Big[\bigsqcup_{g \in A_1 \sqcup A_2}\CClass(g)\Big] &\cong \begin{cases}
        \CSCom(Q_{4n})[\CClass(b)] \vee \CSCom(Q_{4n})[\CClass(ab)], \text{ if } n \text{ is odd} \\
        \CSCom(Q_{4n})[\CClass(b)] \sqcup \CSCom(Q_{4n})[\CClass(ab)], \text{ if } n \text{ is even}
    \end{cases} \\
    & \cong \begin{cases}
        K_n \vee K_n = K_{2n},  \text{ if } n \text{ is odd }\\
K_n \sqcup K_n=2K_n,  \text{ if } n \text{ is even.}
    \end{cases}
\end{align*}}
Thus by Lemma \ref{Lemma1} and after the above discussion, we have
\begin{align*}
    \CSCom(Q_{4n})  \cong K_2 \vee \CSCom(Q_{4n})\Big[ \bigsqcup_{g \in A}\CClass(g) \Big]
    & \cong K_2 \vee \CSCom(Q_{4n})\Big[ \bigsqcup_{g \in A_1 \sqcup A_2 \sqcup A_3}\CClass(g) \Big] \\
    & \cong \begin{cases}
        K_2 \vee \left(K_{2n-2} \sqcup 2K_n\right), & \text{ if } n  \text{ is even} \\
            K_2 \vee \left(K_{2n-2} \sqcup K_{2n} \right), & \text{ if } n \text{ is odd.}
    \end{cases}
\end{align*}
\end{proof}
\begin{theorem}\label{CSCom V_8n}
     If $V_{8n} = \langle a, b: a^{2n}=b^{4}=1, ba=a^{-1}b^{-1}, b^{-1}a=a^{-1}b\rangle$ then 
\begin{align*}
    \CSCom(V_{8n}) \cong \begin{cases}
         K_4 \vee \left(K_{4n-4} \sqcup 2K_{2n}\right), & \text{ if } n \text{ is even} \\
         K_2 \vee \left(K_{4n-2} \sqcup 2K_{2n}\right), & \text{ if } n \text{ is odd.}
    \end{cases} 
\end{align*}
\end{theorem}
\begin{proof}
If $n$ is even then $Z(V_{8n})=\{1, a^n, b^2, a^nb^2\}$, $\CClass(b)=\left\{a^{4t}b, a^{4t+2}b^3 | 1 \leq t \leq \frac{n}{2}\right\}$, $\CClass(ab)=\left\{a^{4t+1}b, a^{4t+3}b^3 \, | \, 1 \leq t \leq \frac{n}{2}\right\}$, \, \, $\CClass(b^3)=\left\{a^{4t+2}b, a^{4t}b^3|1 \leq t \leq \frac{n}{2}\right\}$, $\CClass(ab^3)=\left\{a^{4t+3}b, a^{4t+1}b^3|\right.$ $\left.1 \leq t \leq \frac{n}{2}\right\}$, $\CClass(a^{2l})=\left\{a^{2l}, a^{-2l}\right\}$, $\CClass(a^{2m-1})=\left\{a^{2m-1}, a^{-2m+1}b^2\right\}$ and $\CClass(a^{2s}b^2)=\left\{a^{2s}b^2, a^{-2s}b^2\right\}$, where $1 \leq l \leq \frac{n}{2}-1, 1 \leq m \leq n$ and $1 \leq s \leq \frac{n}{2}-1$, are the distinct conjugacy classes of $V_{8n}$. As such, we have $\CSCom(V_{8n})[\CClass(b)] \cong K_n$, $\CSCom(V_{8n})[\CClass(ab)] \cong K_n$, $\CSCom(V_{8n})[\CClass(b^3)] \cong K_n$ and $\CSCom(V_{8n})[\CClass(ab^3)] \cong K_n$. Let $A=\{b, ab, b^3, ab^3, a^2, a^4, \ldots, a^{n-2}, a, a^3, a^5, \ldots, a^{2n-1}, a^2b^2, a^4b^2, \ldots, a^{n-2}b^2\}$, $A_1=\{b, b^3\}$, \, $A_2=\{ab, ab^3\}$ \, and \,  $A_3=\{a^2, a^4, \ldots, a^{n-2}, a, a^3, a^5, \ldots, a^{2n-1}, a^2b^2,$ $ a^4b^2, \ldots,$ \, $a^{n-2}b^2\}$. Since all powers of $a$ commute with each other and $b^2 \in Z(V_{8n})$, therefore, $\CSCom(V_{8n})[\bigsqcup\limits_{g\in A_3} \CClass(g)] \cong \underbrace{K_2 \vee \cdots \vee K_2}_{(\frac{n}{2}-1+n+\frac{n}{2}-1)\text{-times}}=K_{4n-4}$. For all $x \in A_1 \sqcup A_2$ and $y \in A_3$ and for all $x' \in \CClass(x)$ and $y' \in \CClass(y)$ we have $x'y' \neq y'x'$. Also, since $b$ commutes with $b^3$ and $ab$ commutes with $ab^3$, therefore, $\CSCom(V_{8n})[\bigsqcup\limits_{g\in A_1} \CClass(g)] \cong K_n \vee K_n=K_{2n}$ and $\CSCom(V_{8n})[\bigsqcup\limits_{g\in A_2} \CClass(g)] \cong K_n \vee K_n=K_{2n}$ respectively. However, for all $x \in A_1$ and $y \in A_2$ for all $x' \in \CClass(x)$ and $y' \in \CClass(y)$ we have $x'y' \neq y'x'$. Therefore, $\CSCom(V_{8n})[\bigsqcup\limits_{g\in A_1 \sqcup A_2} \CClass(g)] \cong K_{2n} \sqcup K_{2n} = 2K_{2n}$.

 Thus by Lemma \ref{Lemma1} and after the above discussion, we have
\begin{align*}
    \CSCom(V_{8n})  \cong K_4 \vee \CSCom(V_{8n}) \Big[ \bigsqcup_{g \in A}\CClass(g) \Big]
    & \cong K_4 \vee \CSCom(V_{8n})\Big[ \bigsqcup_{g \in A_1 \sqcup A_2 \sqcup A_3}\CClass(g) \Big] \\
    & \cong K_4 \vee (2K_{2n} \sqcup K_{4n-4}).
\end{align*}

If $n$ is odd then $Z(V_{8n})=\{1, b^2\}$, $\CClass(b)=\left\{a^{2t}b, a^{2t}b^3 | 1 \leq t \leq n\right\}$, $\CClass(ab)=\left\{a^{2t+1}b, \right.$ $\left. a^{2t+1}b^3 | 1 \leq t \leq n\right\}$, $\CClass(a^{2l})=\left\{a^{2l}, a^{-2l}\right\}$, $\CClass(a^{2m-1})=\left\{a^{2m-1}, a^{-2m+1}b^2\right\}$ and $\CClass(a^{2s}b^2)=\left\{a^{2s}b^2, \right.$ $ \left.a^{-2s}b^2\right\}$, where $1 \leq l \leq \frac{n-1}{2}$, $1 \leq m \leq n$ and $1 \leq s \leq \frac{n-1}{2}$, are the distinct conjugacy classes of $V_{8n}$. As such, we have $\CSCom(V_{8n})[\CClass(b)] \cong K_{2n}$ and $\CSCom(V_{8n})[\CClass(ab)] $ $ \cong K_{2n}$. Let $A=\{b, ab, a^2, a^4, \ldots, a^{n-1}, a, a^3, \ldots, a^{2n-1}, a^2b^2, a^4b^2, \ldots, a^{n-1}b^2\}$, $A_1=\{b\}$, $A_2=\{ab\}$ and $A_3=\{a^2, a^4, \ldots, a^{n-1}, a, a^3, \ldots, a^{2n-1}, a^2b^2, a^4b^2, \ldots, a^{n-1}b^2\}$. Since all powers of $a$ commute with each other and $b^2 \in Z(V_{8n})$, therefore, $\CSCom(V_{8n})[\bigsqcup\limits_{g\in A_3} \CClass(g)] $ $ \cong \underbrace{K_2 \vee \cdots \vee K_2}_{(\frac{n-1}{2}+n+\frac{n-1}{2})\text{-times}}=K_{4n-2}$. For all $x \in A_1 \sqcup A_2$ and $y \in A_3$ and for all $x' \in \CClass(x)$ and $y' \in \CClass(y)$ we have $x'y' \neq y'x'$. Also, for all $x \in A_1$ and $y \in A_2$ and for all $x' \in \CClass(x)$ and $y' \in \CClass(y)$ we have $x'y' \neq y'x'$. Therefore, $\CSCom(V_{8n})[\bigsqcup\limits_{g\in A_1 \sqcup A_2} \CClass(g)] \cong K_{2n} \sqcup K_{2n} = 2K_{2n}$. 

Thus by Lemma \ref{Lemma1} and after the above discussion, we have
\begin{align*}
    \CSCom(V_{8n})  \cong K_2 \vee \CSCom(V_{8n}) \Big[ \bigsqcup_{g \in A}\CClass(g) \Big] 
    & \cong K_2 \vee \CSCom(V_{8n})\Big[ \bigsqcup_{g \in A_1 \sqcup A_2 \sqcup A_3}\CClass(g) \Big] \\
    & \cong K_2 \vee (2K_{2n} \sqcup K_{4n-2}).
\end{align*}
\end{proof}
\begin{theorem}\label{CSCom SD_8n}
Let $SD_{8n} = \langle a, b: a^{4n}=b^{2}=1, bab^{-1}=a^{2n-1} \rangle$ $( n \geq 2)$ be the semi-dihedral group. Then
\begin{align*}
    \CSCom(SD_{8n}) \cong \begin{cases}
         K_2 \vee \left(K_{4n-2} \sqcup 2K_{2n}\right), & \text{ if } n \text{ is even} \\
         K_4 \vee \left(K_{4n-4} \sqcup K_{4n} \right), & \text{ if } n \text{ is odd.}
    \end{cases} 
\end{align*}
\end{theorem}
\begin{proof}
If $n$ is even then $Z(SD_{8n})=\{1, a^{2n}\}$ and $\CClass(b)=\left\{a^{2t}b| 1 \leq t \leq 2n\right\}$, $\CClass(ab)=\left\{a^{2t-1}b|\right.$ $\left. 1 \leq t \leq 2n\right\}$, $\CClass(a^{2k})=\left\{a^{2k}, a^{-2k}\right\}$ and $\CClass(a^{2l+1})=\left\{a^{2l+1}, a^{2n-(2l+1)}\right\}$, where  $1 \leq k \leq n-1$ and $\frac{-n}{2} \leq l \leq \frac{n}{2}-1$, are the distinct conjugacy classes of $SD_{8n}$. As such, we have $\CSCom(SD_{8n})[\CClass(b)] \cong K_{2n}$ and $\CSCom(SD_{8n})[\CClass(ab)] \cong K_{2n}$. Let $A=\{b, ab, a^2, a^4, \ldots, a^{2n-2}, a^{-n+1}, \ldots, a^{n-1}\}$, $A_1=\{b\}$, $A_2=\{ab\}$ and $A_3=\{a^2, a^4, \ldots, a^{2n-2},$ $ a^{-n+1}, \ldots, a^{n-1}\}$. Since all powers of $a$ commute with each other, therefore, $\CSCom(SD_{8n})$ $[\bigsqcup\limits_{g\in A_3} \CClass(g)] $ $ \cong \underbrace{K_2 \vee \cdots \vee K_2}_{(n-1+\frac{n}{2}-1+\frac{n}{2}+1)\text{-times}}=K_{4n-2}$. For all $x \in A_1 \sqcup A_2$ and $y \in A_3$ and for all $x' \in \CClass(x)$ and $y' \in \CClass(y)$ we have $x'y' \neq y'x'$. Also, for all $x \in A_1$ and $y \in A_2$ and for all $x' \in \CClass(x)$ and $y' \in \CClass(y)$ we have $x'y' \neq y'x'$. Therefore, $\CSCom(SD_{8n})[\bigsqcup\limits_{g\in A_1 \sqcup A_2} \CClass(g)] \cong K_{2n} \sqcup K_{2n} = 2K_{2n}$. 

Thus by Lemma \ref{Lemma1} and after the above discussion, we have
\begin{align*}
    \CSCom(SD_{8n})  & \cong K_2 \vee \CSCom(SD_{8n}) \Big[ \bigsqcup_{g \in A}\CClass(g) \Big] \\
    & \cong K_2 \vee \CSCom(SD_{8n})\Big[ \bigsqcup_{g \in A_1 \sqcup A_2 \sqcup A_3}\CClass(g) \Big] \\
    & \cong K_2 \vee (2K_{2n} \sqcup K_{4n-2}).
\end{align*}

If $n$ is odd then $Z(SD_{8n})=\{1, a^n, a^{2n}, a^{3n}\}$ and $\CClass(b)=\left\{a^{4t}b|1 \leq t \right. $ $\left. \leq n \right\}$, $\CClass(ab)=\left\{a^{4t+1}b|1 \leq t \leq n \right\}$, $\CClass(a^2b)=\left\{a^{4t+2}b| 1 \leq t \leq n \right\}$, $\CClass(a^3b)=\left\{a^{4t+3}b|1 \leq t \leq n\right\}$, $\CClass(a^{2k})=\left\{a^{2k}, a^{-2k}\right\}$ and $\CClass(a^{2l+1})=\left\{a^{2l+1}, a^{2n-(2l+1)}\right\}$, where $1 \leq k \leq n-1$ and $-\frac{n-1}{2} \leq l \leq \frac{n-1}{2}-1$, are the distinct conjugacy classes of $SD_{8n}$. As such, we have $\CSCom(SD_{8n})[\CClass(b)] $ $\cong K_n$, $\CSCom(SD_{8n})[\CClass(ab)] \cong K_n$, $\CSCom(SD_{8n})[\CClass(a^2b)] \cong K_n$ and $\CSCom(SD_{8n})$ $[\CClass(a^3b)] \cong K_n$. Let $A=\{b, ab, a^2b, a^3b,  a^2, a^4, \ldots, a^{2n-2}, a^{-n+2}, \ldots, a^{n-2}\}$, $A_1=\{b, ab, $ $a^2b, a^3b\}$ and $A_2=\{a^2, a^4, \ldots, a^{2n-2},$ $ a^{-n+2}, \ldots, a^{n-2}\}$. Since all powers of $a$ commute with each other, therefore, $\CSCom(SD_{8n})$ $[\bigsqcup\limits_{g\in A_2} \CClass(g)] $ $ \cong \underbrace{K_2 \vee \cdots \vee K_2}_{(n-1+\frac{n-1}{2}-1+\frac{n-1}{2}+1)\text{-times}}=K_{4n-4}$. For all $x \in A_1$ and $y \in A_2$ and for all $x' \in \CClass(x)$ and $y' \in \CClass(y)$ we have $x'y' \neq y'x'$. Further, we have $b \in \CClass(b)$ and $a^{2n}b \in \CClass(a^2b)$. Also, if $a^nb \in \CClass(ab)$ then $a^{3n}b \in \CClass(a^3b)$ otherwise if $a^{3n}b \in \CClass(ab)$ then $a^nb \in \CClass(a^3b)$. Since $b, a^nb, a^{2n}b$ and $a^{3n}b$ commute with each other, therefore, $ \CSCom(SD_{8n})[\bigsqcup\limits_{g\in A_1} \CClass(g)] \cong K_n \vee K_n \vee K_n \vee K_n =K_{4n}$.

Thus by Lemma \ref{Lemma1} and after the above discussion, we have
\begin{align*}
    \CSCom(SD_{8n})  \cong K_4 \vee \CSCom(SD_{8n}) \Big[ \bigsqcup_{g \in A}\CClass(g) \Big]
    & \cong K_4 \vee \CSCom(SD_{8n})\Big[ \bigsqcup_{g \in A_1 \sqcup A_2}\CClass(g) \Big] \\
    & \cong K_4 \vee (K_{4n} \sqcup K_{4n-4}).
\end{align*}
\end{proof}
\begin{cor}\label{CSCom QD_2^n}
    Let $QD_{2^m} = \langle a, b : a^{2^{m-1}} = b^2 = 1, bab^{-1} = a^{2^{m-2}-1} \rangle$ $(m\geq 3)$ be the quasi-dihedral group. Then 
    \begin{align*}
        \CSCom(QD_{2^m})\cong K_2 \vee \left( 2K_{2^{m-2}} \sqcup K_{2^{m-1}-2} \right).
    \end{align*}
\end{cor}
\begin{proof}
Notice that $SD_{8n} = QD_{2^m}$ for $n=2^{m-3}$. Therefore, the result follows from Theorem \ref{CSCom SD_8n}.
\end{proof}

\begin{theorem}\label{CSCom U_6n}
    If $U_{6n}= \langle a, b: a^{2n}=b^3=1, a^{-1}ba=b^{-1} \rangle$ then
    \[
    \CSCom(U_{6n}) \cong K_n \vee (K_{2n} \sqcup K_{3n}).
    \]
\end{theorem}
\begin{proof}
We know, $Z(U_{6n})=\{1, a^2, a^4, \ldots, a^{2n-2}\}$ and $\CClass(a^{2m-1})=\left\{a^{2m-1}, a^{2m-1}b, a^{2m-1}b^2\right\}$ and $\CClass(a^{2s}b)=\left\{a^{2s}b, a^{2s}b^2\right\}$, where $1 \leq m \leq n$ and $1 \leq s \leq n$, are the distinct conjugacy classes of $U_{6n}$. Let $A=\{a, a^3, \ldots, a^{2n-1}, b, a^2b, a^4b, \ldots, a^{2n-2}b\}$, $A_1=\{a, a^3, \ldots, a^{2n-1}\}$ and $A_2=\{b, a^2b, a^4b, \ldots, a^{2n-2}b\}$. Since all powers of $a$ commute with each other, therefore, $\CSCom(U_{6n})$ $[\bigsqcup\limits_{g\in A_1} \CClass(g)] $ $ \cong \underbrace{K_3 \vee \cdots \vee K_3}_{n\text{-times}}=K_{3n}$. Also, since $a^{2l} \in Z(U_{6n})$ for $1 \leq l \leq n$, therefore, $\CSCom(U_{6n})$ $[\bigsqcup\limits_{g\in A_2} \CClass(g)] $ $ \cong \underbrace{K_2 \vee \cdots \vee K_2}_{n\text{-times}}=K_{2n}$. However, for all $x \in A_1$ and $y \in A_2$ and for all $x' \in \CClass(x)$ and $y' \in \CClass(y)$ we have $x'y' \neq y'x'$.

Thus by Lemma \ref{Lemma1} and after the above discussion, we have
\begin{align*}
    \CSCom(U_{6n})  \cong K_{n} \vee \CSCom(U_{6n}) \Big[ \bigsqcup_{g \in A}\CClass(g) \Big]
    & \cong K_n \vee \CSCom(U_{6n})\Big[ \bigsqcup_{g \in A_1 \sqcup A_2}\CClass(g) \Big] \\
    & \cong K_n \vee (K_{3n} \sqcup K_{2n}).
\end{align*}
\end{proof}
\begin{theorem}\label{CSCom M_2mn}
Let $M_{2mn}= \langle a, b: a^{m}=b^{2n}=1, bab^{-1}=a^{-1} \rangle$ ($m \geq 3$ but not equal to $4$) be a metacyclic group. Then
\begin{align*}
    \CSCom(M_{2mn}) \cong \begin{cases}
        K_n \vee \left(K_{mn-n} \sqcup K_{mn}\right), & \text{ if } m \text{ is odd} \\
        K_{2n} \vee \left(2K_{\frac{mn}{2}} \sqcup K_{mn-2n}\right), & \text{ if } m \text{ and } \frac{m}{2} \text{ are even} \\
        K_{2n} \vee \left(K_{mn} \sqcup K_{mn-2n}\right), & \text{ if } m \text{ is even and } \frac{m}{2} \text{ is odd.}
    \end{cases} 
\end{align*}
\end{theorem}
\begin{proof}
If $m$ is odd then $Z(M_{2mn})=\{1, b^2, b^4, \ldots, b^{2n-2}\}$ and $\CClass(a^s)=\left\{a^s, a^{m-s}\right\}$, $\CClass(b^{2u-1})$ $=\left\{a^tb^{2u-1}| 1 \leq t \leq \right.$ $\left. m \right\}$ \quad and $\CClass(ab^{2r}) =$ $\left\{a^tb^{2r}| 1 \leq t \leq m-1\right\}$,  where $1 \leq r \leq n-1, 1 \leq s \leq \frac{m-1}{2}$ and $1 \leq u \leq n$, are the distinct conjugacy classes of $M_{2mn}$. Let $A=\{a, a^2, \ldots, a^{\frac{m-1}{2}}, b, b^3, \ldots, $ $ b^{2n-1}, ab^2, ab^4, \ldots, ab^{2n-2}\}$, \, $A_1=\{a, a^2, \ldots, a^{\frac{m-1}{2}}, ab^2, ab^4, \ldots, $ $ab^{2n-2}\}$ \, and \, $A_2=\{b, b^3, \ldots, b^{2n-1}\}$. Since all powers of $a$ commute with each other and $b^{2r} \in Z(M_{2mn})$ for $1 \leq r \leq n-1$, therefore, \begin{align*}
\CSCom(M_{2mn})[\bigsqcup\limits_{g\in A_1} \CClass(g)] \cong \underbrace{K_2 \vee \cdots \vee K_2}_{\frac{m-1}{2}\text{-times}} \vee \underbrace{K_{m-1} \vee \cdots \vee K_{m-1}}_{(n-1)\text{-times}}&=K_{m-1} \vee K_{mn-m-n+1} \\
&=K_{mn-n}.
\end{align*}
Since all powers of $b$ commute with each other, therefore, $\CSCom(M_{2mn})$ $[\bigsqcup\limits_{g\in A_2} \CClass(g)] \cong \underbrace{K_m \vee \cdots \vee K_m}_{n\text{-times}}=K_{mn}$. However, for all $x \in A_1$ and $y \in A_2$ and for all $x' \in \CClass(x)$ and $y' \in \CClass(y)$ we have $x'y' \neq y'x'$. 

Thus by Lemma \ref{Lemma1} and after the above discussion, we have
\begin{align*}
    \CSCom(M_{2mn}) &\cong K_{n} \vee \CSCom(M_{2mn}) \Big[ \bigsqcup_{g \in A}\CClass(g) \Big] \\
    & \cong K_n \vee \CSCom(M_{2mn})\Big[\bigsqcup_{g \in A_1 \sqcup A_2}\CClass(g) \Big] \\
    & \cong K_n \vee (K_{mn-n} \sqcup K_{mn}).
\end{align*}
     
If $m$ is even and $\frac{m}{2}$ is odd then $Z(M_{2mn})=\{1, b^2, b^4, \ldots, b^{2n-2}, a^{\frac{m}{2}}, a^{\frac{m}{2}}b^2, a^{\frac{m}{2}}b^4, \ldots,$ $ a^{\frac{m}{2}}b^{2n-2}\}$ and $\CClass(a^h)=\left\{a^h, a^{m-h}\right\}$, $\CClass(b^{2u-1})=\left\{a^{2t}b^{2u-1}|1 \leq t \leq \frac{m}{2}\right\}$, $\CClass(ab^{2u-1})=\left\{a^{2t-1}b^{2u-1}|1 \leq t \leq \frac{m}{2}\right\}$, \, $\CClass(ab^{2s})=\left\{a^{2t-1}b^{2s}|1 \leq t \right.$ $\left.\leq \frac{m}{2}, t \neq \frac{m+2}{4}\right\}$ \, and \, $\CClass(a^2b^{2s})=\left\{a^{2t}b^{2s}|1 \leq t \leq \frac{m}{2}-1\right\}$, where $1 \leq s \leq n-1, 1 \leq u \leq n$ and $1 \leq h \leq \frac{m}{2}-1$, are the distinct conjugacy classes of $M_{2mn}$. Let $A=\{a, a^2, \ldots, a^{\frac{m}{2}-1}, b, b^3, \ldots, b^{2n-1}, ab, ab^3, \ldots,$ $ ab^{2n-1}, ab^2, ab^4, \ldots, ab^{2n-2}, a^2b^2, a^2b^4, \ldots, a^2b^{2n-2}\}$, \, $A_1=\{a, a^2, \ldots, a^{\frac{m}{2}-1}, ab^2, ab^4, \ldots,$  $ ab^{2n-2}, a^2b^2, a^2b^4, \ldots, a^2b^{2n-2}\}$ and $A_2=\{b, b^3, \ldots, b^{2n-1}, ab, ab^3, \ldots, ab^{2n-1}\}$. Since all powers of $a$ commute with each other and $b^{2s} \in Z(M_{2mn})$ for $1 \leq s \leq n-1$, therefore, 
\begin{align*}
\CSCom(M_{2mn})[\bigsqcup\limits_{g\in A_1} \CClass(g)] & \cong \underbrace{K_2 \vee \cdots \vee K_2}_{(\frac{m}{2}-1)\text{-times}} \vee \underbrace{K_{\frac{m}{2}-1} \vee \cdots \vee K_{\frac{m}{2}-1}}_{(n-1)\text{-times}} \vee \underbrace{K_{\frac{m}{2}-1} \vee \cdots \vee K_{\frac{m}{2}-1}}_{(n-1)\text{-times}} \\
&=K_{m-2} \vee K_{\frac{mn}{2}-\frac{m}{2}-n+1}  \vee K_{\frac{mn}{2}-\frac{m}{2}-n+1} = K_{mn-2n}.
\end{align*}
Now, let $x \in \bigsqcup\limits_{u=1}^{n}\CClass(b^{2u-1})$ and $y \in \bigsqcup\limits_{u=1}^{n}\CClass(ab^{2u-1})$ then $x$ is adjacent to $y$ for all $x$ and $y$ since $a^{2t}b^{2u-1} \in \CClass(x)$ commutes with $a^{2t-1}b^{2u-1} \in \CClass(y)$ for some $t \in \{1, \ldots, \frac{m}{2}\}$ and $u \in \{1, \ldots, n\}$. Also, all powers of $b$ commute with each other and $ab^{2u-1}$ commutes with $ab^{2u'-1}$ for all $u, u' \in \{1, 2, \ldots, n\}$. Therefore
\begin{align*}
\CSCom(M_{2mn})[\bigsqcup\limits_{g\in A_2} \CClass(g)]  \cong \underbrace{K_{\frac{m}{2}} \vee \cdots \vee K_{\frac{m}{2}}}_{n\text{-times}} \vee \underbrace{K_{\frac{m}{2}} \vee \cdots \vee K_{\frac{m}{2}}}_{n\text{-times}} 
&= K_{\frac{mn}{2}+\frac{mn}{2}} \\
&=K_{mn}.
\end{align*}
However, for all $x \in A_1$ and $y \in A_2$ and for all $x' \in \CClass(x)$ and $y' \in \CClass(y)$ we have $x'y' \neq y'x'$. 

Thus by Lemma \ref{Lemma1} and after the above discussion, we have
\begin{align*}
    \CSCom(M_{2mn}) &\cong K_{2n} \vee \CSCom(M_{2mn}) \Big[ \bigsqcup_{g \in A}\CClass(g) \Big] \\
    & \cong K_{2n} \vee \CSCom(M_{2mn})\Big[\bigsqcup_{g \in A_1 \sqcup A_2}\CClass(g) \Big] \\
    & \cong K_{2n} \vee (K_{mn-2n} \sqcup K_{mn}).
\end{align*}
If $m$ and $\frac{m}{2}$ are even then $Z(M_{2mn})=\{1, b^2, b^4, \ldots, b^{2n-2}, a^{\frac{m}{2}}, a^{\frac{m}{2}}b^2, a^{\frac{m}{2}}b^4, \ldots,$ $ a^{\frac{m}{2}}b^{2n-2}\}$ and $\CClass(a^h)=\left\{a^h, a^{m-h}\right\}$, $\CClass(b^{2u-1})=\left\{a^{2t}b^{2u-1}|1 \leq t \leq \frac{m}{2}\right\}$, $\CClass(ab^{2u-1})=\left\{a^{2t-1}b^{2u-1}| \right.$ $\left. 1 \leq t \leq \frac{m}{2}\right\}$, $\CClass(ab^{2s})=\left\{a^{2t-1}b^{2s}|1 \leq t \right.$ $\left.\leq \frac{m}{2}\right\}$ \, and \, $\CClass(a^2b^{2s})=\left\{a^{2t}b^{2s}|1 \leq t \leq \frac{m}{2}-1, \right.$ $\left. t \neq \frac{m}{4}\right\}$, where $1 \leq s \leq n-1, 1 \leq u \leq n$ and $1 \leq h \leq \frac{m}{2}-1$, are the distinct conjugacy classes of $M_{2mn}$. Let $A=\{a, a^2, \ldots, a^{\frac{m}{2}-1}, b, b^3, \ldots, b^{2n-1}, ab, ab^3, \ldots, ab^{2n-1}, ab^2, ab^4, \ldots, $ $ab^{2n-2}, a^2b^2, a^2b^4, \ldots, a^2b^{2n-2}\}$, \, $A_1=\{a, a^2, \ldots, a^{\frac{m}{2}-1}, ab^2, ab^4, \ldots,$  $ ab^{2n-2}, a^2b^2, a^2b^4, \ldots, $ $a^2b^{2n-2}\}$, $A_2=\{b, b^3, \ldots, b^{2n-1}\}$ and $A_3=\{ab, ab^3, \ldots, ab^{2n-1}\}$. Since all powers of $a$ commute with each other and $b^{2s} \in Z(M_{2mn})$ for $1 \leq s \leq n-1$, therefore, 
\begin{align*}
\CSCom(M_{2mn})[\bigsqcup\limits_{g\in A_1} \CClass(g)] & \cong \underbrace{K_2 \vee \cdots \vee K_2}_{(\frac{m}{2}-1)\text{-times}} \vee \underbrace{K_{\frac{m}{2}} \vee \cdots \vee K_{\frac{m}{2}}}_{(n-1)\text{-times}} \vee \underbrace{K_{\frac{m}{2}-2} \vee \cdots \vee K_{\frac{m}{2}-2}}_{(n-1)\text{-times}} \\
&=K_{m-2} \vee K_{\frac{mn}{2}-\frac{m}{2}}  \vee K_{\frac{mn}{2}-\frac{m}{2}-2n+2} = K_{mn-2n}.
\end{align*}
Since all powers of $b$ commute with each other, we have $\CSCom(M_{2mn})[\bigsqcup\limits_{g\in A_2} \CClass(g)] \cong \underbrace{K_{\frac{m}{2}} \vee \cdots \vee K_{\frac{m}{2}}}_{n\text{-times}}=K_{\frac{mn}{2}}$. Also $ab^{2u-1}$ commutes with $ab^{2u'-1}$ for all $u, u' \in \{1, 2, \ldots, n\}$ so we have $\CSCom(M_{2mn})[\bigsqcup\limits_{g\in A_3} \CClass(g)] \cong \underbrace{K_{\frac{m}{2}} \vee \cdots \vee K_{\frac{m}{2}}}_{n\text{-times}}=K_{\frac{mn}{2}}$. However, for all $x \in A_2$ and $y \in A_3$ and for all $x' \in \CClass(x)$ and $y' \in \CClass(y)$ we have $x'y' \neq y'x'$. Similarly, for all $x \in A_1$ and $y \in A_2 \sqcup A_3$ and for all $x' \in \CClass(x)$ and $y' \in \CClass(y)$ we have $x'y' \neq y'x'$. Thus by Lemma \ref{Lemma1} and after the above discussion, we have
\begin{align*}
    \CSCom(M_{2mn}) &\cong K_{2n} \vee \CSCom(M_{2mn}) \Big[ \bigsqcup_{g \in A}\CClass(g) \Big] \\
    & \cong K_{2n} \vee \CSCom(M_{2mn})\Big[\bigsqcup_{g \in A_1 \sqcup A_2 \sqcup A_3}\CClass(g) \Big] \\
    & \cong K_{2n} \vee (K_{mn-2n} \sqcup K_{\frac{mn}{2}} \sqcup K_{\frac{mn}{2}})= K_{2n} \vee (K_{mn-2n} \sqcup 2K_{\frac{mn}{2}}).
\end{align*}
\end{proof}
\subsection{Order super commuting graphs}
The order super commuting graph of a group $G$, denoted by $\OSCom(G)$, is a graph with vertex set $v(\OSCom(G))=G$ and two distinct vertices $g$ and $h$ are adjacent if either $[g]=[h]$ or there exist $g' \in [g]$ and $h' \in [h]$ such that $g'h'=h'g'$, where $[g] = \OrdClass(g)$ is the order class of $g$ in $G$ (that is, $[g]$ contains all elements having same order as $g$ in $G$). In this section, we obtain the graph structures of $\OSCom(G)$
 for the groups considered in Subsection \ref{CScom}.
\begin{theorem}\label{OSCom D_2n}
Let $D_{2n}=\langle a, b: a^n=b^2=1, bab^{-1}=a^{-1} \rangle$ $(n \geq 3)$ be the dihedral group. Then
\begin{align*}
\OSCom(D_{2n}) \cong \begin{cases}
           K_{2n}, & \text{ if } n \text{ is even} \\
           K_1 \vee (K_{n-1} \sqcup K_n), & \text{ if } n \text{ is odd.} 
       \end{cases}
\end{align*}
\end{theorem}
\begin{proof}
Let $x$ be an arbitrary element of $D_{2n}$. If $x=1$ then clearly $x \sim y$ for all $1 \neq y \in D_{2n}$.
If $n$ is even then $X_1= \{ a^{\frac{n}{2}}, b, ab, \ldots, a^{n-1}b\}$ forms an order class having elements of order $2$. 
If $x \in X_1$ then $a^{\frac{n}{2}} \in \OrdClass(x)$ which commute with all other elements of $D_{2n}$. Therefore, $x \sim y$ for all $y \neq x \in D_{2n}$.
If $x \in X_2 = \{a, a^2, \ldots, a^{n-1}\}\setminus \{a^{\frac{n}{2}}\}$ then $x \sim y$ for all $y \neq x \in D_{2n}$ since $a^j \in \OrdClass(x)$ for some $j$ and $a^j$ commutes with $a^{\frac{n}{2}}$ and all other powers of $a$. 
Hence
\[ 
\OSCom(D_{2n}) \cong K_{2n}.
\]

If $n$ is odd then $Y_1=\{b, ab, \ldots, a^{n-1}b\}$ forms an order class having elements of order $2$. As such, $\OSCom(D_{2n})[Y_1] \cong K_n$. If $x \in Y_2=\{a, a^2, \ldots, a^{n-1}\}$ then $x \sim y$ for all $x \neq y \in Y_2$ since all powers of $a$ commute with each other but $x \nsim y$ for all $y \in Y_1$. Therefore,
 $\OSCom(D_{2n})[Y_2] \cong K_{n-1}$ and so $\OSCom(D_{2n})[Y_1 \sqcup Y_2] \cong K_{n-1} \sqcup K_n$.
Hence, the result follows.
\end{proof}
\begin{theorem}\label{OSCom Q_4n}
    Let $Q_{4n}=\langle a, b: a^{2n}=1, a^n=b^2, bab^{-1}=a^{-1} \rangle$ $(n \geq 2)$ be the generalized quaternion group. Then
    \begin{align*}
        \OSCom(Q_{4n}) \cong \begin{cases}
            K_{4n}, & \text{ if } n \text{ is even} \\
             K_2 \vee (K_{2n} \sqcup K_{2n-2}), & \text{ if } n \text{ is odd.} 
             \end{cases}    
    \end{align*}
\end{theorem}
\begin{proof}
      Let $x$ be an arbitrary element of $Q_{4n}$. If $x \in \{1, a^n\}=Z(Q_{4n})$ then clearly $x \sim y$ for all $x \neq y \in Q_{4n}$.
      
      If $n$ is even then $X_1=\{a^{\frac{n}{2}}, b, ab, a^2b, \ldots, a^{2n-1}b\}$ forms an order class having elements of order 4. Let $X_2=\{a, a^2, a^3, \ldots, a^{2n-1}\}\setminus \{a^{\frac{n}{2}}, a^n\}$. For all $x \in X_1$ and $y \in X_2$, we have $x \sim y$ since $a^{\frac{n}{2}} \in \OrdClass(x)$ commutes with $a^j \in \OrdClass(y)$ for all $j$. Also, if $x \in X_2$ then $x \sim y$ for all $x \neq y \in X_2$ since $a^j \in \OrdClass(x)$ for some $j$ commutes with all other powers of $a$. Therefore
     \[
     \OSCom(Q_{4n}) \cong K_{4n}.
     \]
    
     If $n$ is odd then $Y_1=\{b, ab, a^2b, \ldots, a^{2n-1}b\}$ forms an order class having elements of order 4. As such, $\OSCom(Q_{4n})[Y_1] \cong K_{2n}$. Let $ Y_2 =\{a, a^2, a^3, \ldots, a^{2n-1}\}\setminus \{a^n\}$. If $x \in Y_2$ then $x \sim y$ for all $x\neq y \in Y_2$ since all powers of $a$ commute with each other but $x \nsim y$ for all $y \in Y_1$.  Therefore $\OSCom(Q_{4n})[Y_2] \cong K_{2n-2}$
    and so $\OSCom(Q_{4n})[Y_1 \sqcup Y_2] \cong K_{2n} \sqcup K_{2n-2}$. Hence, the result follows.
\end{proof}
\begin{theorem}\label{OSCom V_8n}
     If $V_{8n} = \langle a, b: a^{2n}=b^{4}=1, ba=a^{-1}b^{-1}, b^{-1}a=a^{-1}b\rangle$ then 
     \begin{align*}
     \OSCom(V_{8n}) \cong \begin{cases}
          K_{8n}, & \text{ if } n \text{ is even} \\
          K_{2n+4} \vee (K_{2n} \sqcup K_{4n-4}), & \text{ if } n \text{ is odd.} 
     \end{cases}
     \end{align*}
     \end{theorem}
\begin{proof}
Let $x$ be any arbitrary element of $V_{8n}$. If $x=1$ then $x \sim y$ for all $1 \neq y \in V_{8n}$. 
If $n$ is even, we have
    \[
    X_1=\left\{b^2, a^n, a^nb^2, ab, a^3b, a^5b, \ldots, a^{2n-1}b, ab^3, a^3b^3, a^5b^3, \ldots, a^{2n-1}b^3\right\}
    \]
    and
    \[
    X_2= \left\{a^{\frac{n}{2}}, b, b^3, a^2b, a^4b, \ldots, a^{2n-2}b, a^2b^3, a^4b^3, \ldots, a^{2n-2}b^3 \right\}
    \]
    forming order classes having elements of order 2 and 4 respectively. If $x \in X_1$ then $x \sim y$ for all $x \neq y \in V_{8n}$ since $b^2 \in \OrdClass(x)$ commutes with all other elements of $V_{8n}$.  If $x \in X_2$ then $x \sim y$ for all $y \in X_3=\left\{a, a^2, a^3, \ldots, a^{2n-1}, ab^2, a^2b^2, \right. $ $\left.a^3b^2, \ldots, a^{2n-1}b^2 \right\} \setminus \{a^n, a^nb^2, a^{\frac{n}{2}}\}$ since $a^\frac{n}{2} \in \OrdClass(x)$ commutes with all elements of $X_3$. If $x \in X_3$ then $x \sim y$ for all $x \neq y \in X_3$ since all powers of $a$ commute with each other and $b^2 \in Z(V_{8n})$. Therefore,
      $
      \OSCom(V_{8n}) \cong K_{8n}.
      $
      
    If $n$ is odd, we have
    \[
    Y_1=\left\{b^2, a^n, a^nb^2, ab, a^3b, a^5b, \ldots, a^{2n-1}b, ab^3, a^3b^3, a^5b^3, \ldots, a^{2n-1}b^3\right\}
    \]
    and
    \[
    Y_2= \left\{b, b^3, a^2b, a^4b, \ldots, a^{2n-2}b, a^2b^3, a^4b^3, \ldots, a^{2n-2}b^3 \right\}
    \]
    forming order classes having elements of order 2 and 4 respectively. If $x \in Y_1$ then $x \sim y$ for all $x \neq y \in V_{8n}$ since $b^2 \in \OrdClass(x)$ commutes with all other elements of $V_{8n}$. That is, $Y_1 \sqcup \{1\}$ is a set of dominant vertices in $\OSCom(V_{8n})$ and $\OSCom(V_{8n})[Y_1 \sqcup \{1\}] \cong K_{2n+4}$.
    
     Let $Y_3=\left\{a, a^2, a^3, \ldots, a^{2n-1}, ab^2, a^2b^2, a^3b^2,\right.$ $\left. \ldots, a^{2n-1}b^2 \right\} \setminus \{a^n, a^nb^2\}$. If $x \in Y_2$ then $x \nsim y$ for all $y \in Y_3$. If $x \in Y_3$ then $x \sim y$ for all $x \neq y \in Y_3$ since since all powers of $a$ commute with each other and $b^2 \in Z(V_{8n})$. Therefore, $\OSCom(V_{8n})[Y_2] \cong K_{2n}$, $ \OSCom(V_{8n})[Y_3] \cong K_{4n-4}$
       and 
       \[
       \OSCom(V_{8n})[Y_2 \sqcup Y_3] \cong K_{2n} \sqcup K_{4n-4}.
       \]
       Hence, the result follows.
\end{proof}
\begin{theorem}\label{OSCom U_6n}
    If $U_{6n}= \langle a, b: a^{2n}=b^3=1, a^{-1}ba=b^{-1} \rangle$ then
    \[
    \OSCom(U_{6n}) \cong K_n \vee (K_{2n} \sqcup K_{3n}).
    \]
\end{theorem}
\begin{proof}
    Let $x$ be an arbitrary element of $U_{6n}$. If $x \in \left\{ 1, a^2, a^4, \ldots, a^{2n-2} \right\}=Z(U_{6n})$ then clearly $x \sim y$ for all $x \neq y \in U_{6n}$. That is, $Z(U_{6n})$ is a set of dominant vertices in $\OSCom(U_{6n})$ and $\OSCom(U_{6n})[Z(U_{6n})] \cong K_n$. We have $X_1= \{ a, ab, ab^2 \}$ forming an order class having elements of order $2n$ and $X_i=\left\{ a^{2i-1}, a^{2i-1}b, a^{2i-1}b^2\right\}$ forming order classes where $2 \leq i \leq n$. 

     If $x \in \bigsqcup\limits_{i=1}^{n} X_i$ then $x \sim y$ for all $x \neq y \in \bigsqcup_{i=1}^{n} X_i$ since $a^{2i-1} \in \OrdClass(x)$ commutes with $a^{2i'-1} \in \OrdClass(y)$ for some $i, i' \in \left\{1, \ldots, n\right\}$. Also, $x \nsim y$ for all $y \in X_{n+1}=\left \{a^2b, a^4b, \ldots, a^{2n}b, a^2b^2, a^4b^2,\right.$ $\left. \ldots, a^{2n}b^2 \right\}$. If $x \in X_{n+1}=\left\{a^2b, a^4b, \ldots, a^{2n}b, a^2b^2, a^4b^2, \ldots,\right.$ $\left. a^{2n}b^2 \right\}$ then $x \sim y$ for all $x \neq y \in X_{n+1}$ since $a^{2s} \in Z(U_{6n})$ for $1 \leq s \leq n$. Therefore,
       $
       \OSCom(U_{6n})[\bigsqcup\limits_{i=1}^{n} X_i] \cong \underbrace{K_3 \vee \cdots \vee K_3}_{n\text{-times}}=K_{3n}$ and $\OSCom(U_{6n})[X_{n+1}] \cong K_{2n}.
       $
       Hence, the result follows.
\end{proof}
\begin{theorem}\label{OSCom SD_8n}
Let $SD_{8n} = \langle a, b: a^{4n}=b^{2}=1, bab^{-1}=a^{2n-1} \rangle$  $(n \geq 2)$ be the semi-dihedral group. Then
\[
\OSCom(SD_{8n}) \cong K_{8n}.
\]
\end{theorem}
\begin{proof}
    Let $x$ be an arbitrary element of $SD_{8n}$. If $x=1$ then $x \sim y$ for all $1 \neq y \in SD_{8n}$. We have 
    \[
    X_1=\left\{b, a^{2n}, a^2b, a^4b, \ldots, a^nb, \ldots, a^{2n}b, \ldots, a^{4n-2}b \right\}
    \]
    and 
    \[
    X_2=\left\{a^n, a^{3n}, ab, a^3b, a^5b, \ldots, a^{2n-1}b, a^{2n+1}b, \ldots, a^{3n}b, \ldots, a^{4n-1}b \right\}
    \]
    forming order classes having elements of order 2 and 4 respectively. If $x \in X_1$ then $x \sim y$ for all $x \neq y \in SD_{8n}$ since $a^{2n} \in \OrdClass(x)$ commutes with all other elements of $SD_{8n}$. Let $X_3=\{a, a^2, a^3, \ldots, a^{4n-1}\} \setminus \{a^n, a^{2n}, a^{3n}\}$. If $x \in X_2$ then $x \sim y$ for all $y \in X_3$ since $a^n \in \OrdClass(x)$ commutes with all elements of $X_3$. If $x \in  X_3=\{a, a^2, a^3, \ldots, a^{4n-1}\} \setminus \{a^n, a^{2n}, a^{3n}\}$ then $x \sim y$ for all $x \neq y \in X_3$ since all powers of $a$ commute with each other. Therefore,
      $
      \OSCom(SD_{8n}) \cong K_{8n}.
      $
\end{proof}
We conclude this section with the following corollary.
\begin{cor}\label{OSCom QD_2^n}
    Let $QD_{2^m}=\langle a, b: a^{2^{m-1}}=b^2=1, bab^{-1}=a^{2^{m-2}-1} \rangle$ $(n \geq 4)$ be the quasi-dihedral group. Then
    \[
    \OSCom(QD_{2^m}) \cong K_{2^m}.
    \]
\end{cor}
\section{Zagreb indices of super graphs}
In this section we shall compute the Zagreb indices of various super graphs obtained in Section 2  and check the validity of Hansen--Vuki{\v{c}}evi{\'c} conjecture. The following lemma can be easily verified and is useful in this context.
\begin{lemma}\label{Zagreb}
If $\Gamma=K_a \vee (dK_b \sqcup K_c)$, where $a,b,c,d$ are positive integers, then 
\[
M_1(\Gamma)=a(a+db+c-1)^2+db(a+b-1)^2+c(a+c-1)^2
\]
and
\begin{align*}
	M_2(\Gamma)&=\frac{1}{2}a(a-1)(a+db+c-1)^2+dab(a+db+c-1)(a+b-1) \\
	&  +ac(a+db+c-1)(a+c-1)+\frac{1}{2}db(b-1)(a+b-1)^2+\frac{1}{2}c(c-1)(a+c-1)^2.
\end{align*}
\end{lemma}
\begin{theorem}
The graphs $\ESCom(D_{2n})$, $\CSCom(D_{2n})$ and $\OSCom(D_{2n})$ satisfy the Hansen--Vuki{\v{c}}evi{\'c} conjecture \eqref{Conj-eq}.
\end{theorem}
\begin{proof}
{\it Part 1.} Suppose that $\Gamma(D_{2n})= \ESCom(D_{2n})$. If $n$ is odd, then, by Proposition \ref{ESCom(G)}(a), $\Gamma(D_{2n})\cong K_1 \vee (nK_1 \sqcup K_{n-1})$. Now, by Lemma \ref{Zagreb}, we have 
\[
M_1(\Gamma(D_{2n})) = n^3+n^2 \text{ and }  M_2(\Gamma(D_{2n})) = \frac{n^4-n^3+3n^2-n}{2}.
\]
Also, $|v(\Gamma(D_{2n}))|=|D_{2n}|=2n$ and $|e(\Gamma(D_{2n}))|=n+(n-1)+\binom{n-1}{2}=\frac{n^2+n}{2}$. Therefore
\begin{align*}
	\frac{M_2(\Gamma(D_{2n}))}{|e(\Gamma(D_{2n}))|}- \frac{M_1(\Gamma(D_{2n}))}{|v(\Gamma(D_{2n}))|}&=\frac{n^5-4n^4+5n^3-2n^2}{2n(n^2+n)} \\
	&= \frac{n^4(n-4)+n^2(5n-2)}{2n(n^2+n)}.
\end{align*}
It is easy to see that $n^4(n-4)+n^2(5n-2)>0$ and $2n(n^2+n)>0$ for all $n \geq 3$ such that $n$ is odd. Thus the the Hansen--Vuki{\v{c}}evi{\'c} conjecture holds.

 If $n$ is even, then, by Proposition \ref{ESCom(G)}(a), $\Gamma(D_{2n}) \cong K_2 \vee (\frac{n}{2}K_2 \sqcup K_{n-2})$. Now, by Lemma \ref{Zagreb}, we have 
\[
M_1(\Gamma(D_{2n})) = n^3+4n^2+6n \text{ and } M_2(\Gamma(D_{2n})) = \frac{n^4+n^3+21n^2}{2}.
\]
Also, $|v(\Gamma(D_{2n}))|=|D_{2n}|=2n$ and $|e(\Gamma(D_{2n}))|=1+\frac{n}{2}+4\frac{n}{2}+2(n-2)+\binom{n-2}{2}=\frac{n^2+4n}{2}$. Therefore
\begin{align*}
	\frac{M_2(\Gamma(D_{2n}))}{|e(\Gamma(D_{2n}))|}- \frac{M_1(\Gamma(D_{2n}))}{|v(\Gamma(D_{2n}))|}&= \frac{n^5-6n^4+20n^3-24n^2}{2n(n^2+4n)} \\
	&= \frac{n^4(n-6)+4n^2(5n-6)}{2n(n^2+4n)}.
\end{align*}
It is easy to see that $n^4(n-6)+4n^2(5n-6)>0$ and $2n(n^2+4n)>0$ for all $n \geq 4$ such that $n$ is even. Thus the the Hansen--Vuki{\v{c}}evi{\'c} conjecture holds.

 {\it Part 2.} Suppose that $\Gamma(D_{2n})= \CSCom(D_{2n})$. If $n$ is odd, then, by Theorem \ref{CSCom D_2n}, we have $\Gamma(D_{2n}) \cong K_1 \vee (K_{n-1} \sqcup K_n)$. Now, by Lemma \ref{Zagreb}, we have 
 \[
 M_1(\Gamma(D_{2n})) = 2n^3+n^2-n \text{ and } M_2(\Gamma(D_{2n})) = \frac{2n^4+2n^3-3n^2+n}{2}.
 \]
Also, $|v(\Gamma(D_{2n}))|=|D_{2n}|=2n$ and $|e(\Gamma(D_{2n}))|=2n-1+\binom{n-1}{2}+ \binom{n}{2}=n^2$. Therefore
\begin{align*}
    \frac{M_2(\Gamma(D_{2n}))}{|e(\Gamma(D_{2n}))|}- \frac{M_1(\Gamma(D_{2n}))}{|v(\Gamma(D_{2n}))|}=\frac{(n-1)^2}{2n} > 0.
\end{align*}
If $n$ and $\frac{n}{2}$ are even, then, by Theorem \ref{CSCom D_2n}, we have $\Gamma(D_{2n}) \cong  K_2 \vee (K_{n-2} \sqcup 2K_{\frac{n}{2}})$. Now, by Lemma \ref{Zagreb}, we have 
\[
M_1(\Gamma(D_{2n})) = \frac{5n^3}{4}+5n^2-2n \text{ and } M_2(\Gamma(D_{2n})) = \frac{9n^4}{16}+\frac{21n^3}{8}+\frac{3n^2}{2}-n.
\]
Also, $|v(\Gamma(D_{2n}))|=|D_{2n}|=2n$ and $|e(\Gamma(D_{2n}))|=1+2(n-2)+\binom{n-2}{2} + 2\times 2 \times  \frac{n}{2}+ \binom{\frac{n}{2}}{2}=\frac{3n^2}{4}+n$. Therefore
\begin{align*}
    \frac{M_2(\Gamma(D_{2n}))}{|e(\Gamma(D_{2n}))|}- \frac{M_1(\Gamma(D_{2n}))}{|v(\Gamma(D_{2n}))|}=\frac{\frac{3n^5}{16}+\frac{n^4}{4}-\frac{n^3}{2}}{\frac{3n^3}{2}+2n^2}=\frac{\frac{n^3}{8}(3n^2+4n-8)}{3n^3+4n^2} > 0.
\end{align*}
If $n$ is even and $\frac{n}{2}$ is odd, then, by Theorem \ref{CSCom D_2n}, we have $\Gamma(D_{2n}) \cong K_2 \vee (K_{n-2} \sqcup K_n)$. Now, by Lemma \ref{Zagreb}, we have 
\[
M_1(\Gamma(D_{2n})) = 2n^3+6n^2-2n \text{ and } M_2(\Gamma(D_{2n})) = n^4+5n^3-n.
\] 
Also, $|v(\Gamma(D_{2n}))|=|D_{2n}|=2n$ and $|e(\Gamma(D_{2n}))|=1+2(n-2)+\binom{n-2}{2}+2n+ \binom{n}{2}=n^2+n$. Therefore
\begin{align*}
    \frac{M_2(\Gamma(D_{2n}))}{|e(\Gamma(D_{2n}))|}- \frac{M_1(\Gamma(D_{2n}))}{|v(\Gamma(D_{2n}))|}&=\frac{n^4+5n^3-n}{n^2+n}-\frac{2n^3+6n^2-2n}{2n}=\frac{2n^3(n-2)}{2n(n^2+n)} > 0.
\end{align*}

{\it Part 3.} Suppose that $\Gamma(D_{2n})= \OSCom(D_{2n})$. If $n$ is even, then, by Theorem \ref{OSCom D_2n}, we have $\Gamma(D_{2n}) \cong  K_{2n}$. As such, $M_1(\Gamma(D_{2n}))=2n(2n-1)^2, M_2(\Gamma(D_{2n}))=\frac{2n(2n-1)^3}{2}, |v(\Gamma(D_{2n}))|=|D_{2n}|=2n$ and $|e(\Gamma(D_{2n}))|=\binom{2n}{2}=n(2n-1)$. Therefore
\[
\frac{M_2(\Gamma(D_{2n}))}{|e(\Gamma(D_{2n}))|}=(2n-1)^2= \frac{M_1(\Gamma(D_{2n}))}{|v(\Gamma(D_{2n}))|}.
\]
If $n$ is odd, then, by Theorem \ref{OSCom D_2n}, we have $\Gamma(D_{2n}) \cong K_1 \vee (K_{n-1} \sqcup K_n)$. As such, the result follows from \it{Part 2}.
\end{proof}
\begin{theorem}
    The graphs $\ESCom(Q_{4n})$, $\CSCom(Q_{4n})$ and $\OSCom(Q_{4n})$ satisfy the Hansen--Vuki{\v{c}}evi{\'c} conjecture \eqref{Conj-eq}.
\end{theorem}
\begin{proof}
    {\it Part 1.} Suppose that $\Gamma(Q_{4n})=\ESCom(Q_{4n})$. By Proposition \ref{ESCom(G)}(b), we have $\Gamma(Q_{4n}) \cong K_2 \vee (n K_2 \sqcup K_{2n-2})$. So, by Lemma \ref{Zagreb}, we have 
    \[
    M_1(\Gamma(Q_{4n})) = 8n^3+16n^2+12n \text{ and } M_2(\Gamma(Q_{4n})) = 8n^4+4n^3+42n^2. 
    \]
    Also, $|v(\Gamma(Q_{4n}))|=|Q_{4n}|=4n$ and $|e(\Gamma(Q_{4n}))|=1+2 \times 2n+n+2(2n-2)+\binom{2n-2}{2}=2n^2+4n$. Therefore, for all $n \geq 2$, we have
    \begin{align*}
        \frac{M_2(\Gamma(Q_{4n}))}{|e(\Gamma(Q_{4n}))|}-\frac{M_1(\Gamma(Q_{4n}))}{|v(\Gamma(Q_{4n}))|}&=\frac{16n^5-48n^4+80n^3-48n^2}{4n(2n^2+4n)} \\
        &=\frac{16n^3(n-1)(n-2)+48n^2(n-1)}{4n(2n^2+4n)} >0.
    \end{align*}

{\it Part 2.} Suppose that $\Gamma(Q_{4n}) = \CSCom(Q_{4n})$. If $n$ is even, then, by Theorem \ref{CSCom Q_4n}, we have $\Gamma(Q_{4n}) \cong K_2 \vee (K_{2n-2} \sqcup 2K_n)$. Now, by Lemma \ref{Zagreb}, we have 
 \[
 M_1(\Gamma(Q_{4n})) = 10n^3+20n^2-4n \text{ and } M_2(\Gamma(Q_{4n})) = 9n^4+21n^3+5n^2-2n.
 \]
    Also, $|v(\Gamma(Q_{4n}))|=|Q_{4n}|=4n$ and $|e(\Gamma(Q_{4n}))|=1+2(2n-2)+\binom{2n-2}{2}+2 \times 2n+2\binom{n}{2}=3n^2+2n$. Therefore, for all $n \geq 2$, we have
    \begin{align*}
        \frac{M_2(\Gamma(Q_{4n}))}{|e(\Gamma(Q_{4n}))|}-\frac{M_1(\Gamma(Q_{4n}))}{|v(\Gamma(Q_{4n}))|}=\frac{6n^5+4n^4-8n^3}{4n(3n^2+2n)} 
        =\frac{6n^5+4n^3(n-2)}{4n(3n^2+2n)} >0.
    \end{align*}
 If $n$ is odd, then, by Theorem \ref{CSCom Q_4n}, we have $\Gamma(Q_{4n}) \cong K_2 \vee (K_{2n-2} \sqcup K_{2n})$. Now, by Lemma \ref{Zagreb}, we have 
 \[
 M_1(\Gamma(Q_{4n})) = 16n^3+16n^2+4n \text{ and } M_2(\Gamma(Q_{4n})) = 16n^4+40n^3-2n.
 \]
    Also, $|v(\Gamma(Q_{4n}))|=|Q_{4n}|=4n$ and $|e(\Gamma(Q_{4n}))|=1+2 \times 2n+\binom{2n}{2}+2(2n-2)+\binom{2n-2}{2}=4n^2+2n$. Therefore, for all $n \geq 2$, we have
    \begin{align*}
        \frac{M_2(\Gamma(Q_{4n}))}{|e(\Gamma(Q_{4n}))|}-\frac{M_1(\Gamma(Q_{4n}))}{|v(\Gamma(Q_{4n}))|}&=\frac{8n^3+20n^2-1}{2n+1}-(4n^2+4n+1)=\frac{8n^2-6n}{2n+1} >0.
    \end{align*}

{\it Part 3.} Suppose that $\Gamma(Q_{4n}) = \OSCom(Q_{4n})$. If $n$ is even, then, by Theorem \ref{OSCom Q_4n}, we have $\Gamma(Q_{4n}) \cong K_{4n}$. As such, $M_1(\Gamma(Q_{4n}))=4n(4n-1)^2$ and $M_2(\Gamma(Q_{4n}))= \binom{4n}{2}(4n-1)^2=\frac{4n(4n-1)^3}{2}$. Also, $|v(\Gamma(Q_{4n}))|=|Q_{4n}|=4n$ and $|e(\Gamma(Q_{4n}))|=\binom{4n}{2}=\frac{4n(4n-1)}{2}$. Therefore
    \begin{align*}
        \frac{M_2(\Gamma(Q_{4n}))}{|e(\Gamma(Q_{4n}))|}=(4n-1)^2=\frac{M_1(\Gamma(Q_{4n}))}{|v(\Gamma(Q_{4n}))|}.
    \end{align*}
If $n$ is odd, then, by Theorem \ref{OSCom Q_4n}, we have $\Gamma(Q_{4n}) \cong K_2 \vee ( K_{2n-2} \sqcup K_{2n})$. As such, the result follows from \it{Part 2}.
\end{proof}
\begin{theorem} \label{Zagreb V_8n}
    The graphs \, $\ESCom(V_{8n})$,  $\CSCom(V_{8n})$ and $\OSCom(V_{8n})$ satisfy the Hansen--Vuki{\v{c}}evi{\'c} conjecture \eqref{Conj-eq}.
\end{theorem}
\begin{proof}
{\it Part 1.} Suppose that $\Gamma(V_{8n})=\ESCom(V_{8n})$. If $n$ is even, then, by Proposition \ref{ESCom(G)}(d), we have $\Gamma(V_{8n}) \cong K_4 \vee (K_{4n-4} \sqcup nK_4)$. Now, by Lemma \ref{Zagreb}, we have 
\[
M_1(\Gamma(V_{8n})) = 64n^3+160n^2+168n \text{ and } M_2(\Gamma(V_{8n})) = 128n^4+160n^3+888n^2+196n.
\]
     Also, $|v(\Gamma(V_{8n}))|=|V_{8n}|=8n$ and $|e(\Gamma(V_{8n}))|=6+4(4n-4)+\binom{4n-4}{2}+16n+6n=8n^2+20n$. Therefore
     \begin{align*}
         \frac{M_2(\Gamma(V_{8n}))}{|e(\Gamma(V_{8n}))|}-\frac{M_1(\Gamma(V_{8n}))}{|v(\Gamma(V_{8n}))|}&= \frac{512n^5-1280n^4+2560n^3-1792n^2}{8n(8n^2+20n)} \\
         &= \frac{n^4(512n-1280)+2560n^3-1792n^2}{8n(8n^2+20n)}.
     \end{align*}
     It is easy to see that $n^4(512n-1280)+2560n^3-1792n^2>0$ and $8n(8n^2+20n)>0$ for all $n \geq 2$ such that $n$ is even.

      If $n$ is odd, then, by Proposition \ref{ESCom(G)}(d), we have $\Gamma(V_{8n}) \cong K_2 \vee (K_{4n-2} \sqcup 2nK_2)$. Now, by Lemma \ref{Zagreb}, we have 
      \[
      M_1(\Gamma(V_{8n})) = 64n^3+64n^2+24n \text{ and } M_2(\Gamma(V_{8n})) = 128n^4+96n^3+104n^2.
      \]
     Also, $|v(\Gamma(V_{8n}))|=|V_{8n}|=8n$ and $|e(\Gamma(V_{8n}))|=1+2(4n-2)+\binom{4n-2}{2}+2n+4 \times 2n=8n^2+8n$. Therefore 
     \begin{align*}
        \frac{M_2(\Gamma(V_{8n}))}{|e(\Gamma(V_{8n}))|}-\frac{M_1(\Gamma(V_{8n}))}{|v(\Gamma(V_{8n}))|}&= \frac{512n^5-256n^4+128n^3-192n^2}{8n(8n^2+8n)} \\
        &= \frac{256n^4(n-1)+n^2(256n^3-192)+128n^3}{8n(8n^2+8n)} >0.
     \end{align*}
 
{\it Part 2.} Suppose that $\Gamma(V_{8n}) = \CSCom(V_{8n})$. If $n$ is even, then, by Theorem \ref{CSCom V_8n}, we have $\Gamma(V_{8n}) \cong K_4 \vee (K_{4n-4} \sqcup 2K_{2n})$. Now, by Lemma \ref{Zagreb}, we have 
\[
M_1(\Gamma(V_{8n})) = 80n^3+208n^2+8n \text{ and } M_2(\Gamma(V_{8n})) = 144n^4+456n^3+356n^2+20n.
\]
      Also, $|v(\Gamma(V_{8n}))|=|V_{8n}|=8n$ and $|e(\Gamma(V_{8n}))|=6+4(4n-4)+\binom{4n-4}{2}+16n+2\binom{2n}{2}=12n^2+12n$. Therefore, for all $n \geq 1$, we have
     \begin{align*}
         \frac{M_2(\Gamma(V_{8n}))}{|e(\Gamma(V_{8n}))|}-\frac{M_1(\Gamma(V_{8n}))}{|v(\Gamma(V_{8n}))|}&=\frac{144n^4+456n^3+356n^2+20n}{12n^2+12n} -\frac{ 80n^3+208n^2+8n}{8n} \\
         &= \frac{192n^5+192n^4+256n^3+64n^2}{8n(12n^2+12n)}>0.
     \end{align*}
If $n$ is odd, then, by Theorem \ref{CSCom V_8n}, we have $\Gamma(V_{8n})\cong K_2 \vee (K_{4n-2} \sqcup 2K_{2n})$. Now, by Lemma \ref{Zagreb}, we have 
\[
M_1(\Gamma(V_{8n})) = 80n^3+80n^2-8n \text{ and } M_2(\Gamma(V_{8n})) = 144n^4+168n^3+20n^2-4n.
\] 
      Also, $|v(\Gamma(V_{8n}))|=|V_{8n}|=8n$ and $|e(\Gamma(V_{8n}))|=1+2(4n-2)+\binom{4n-2}{2}+8n+2n(2n-1)=12n^2+4n$. Therefore, for all $n \geq 1$, we have
     \begin{align*}
        \frac{M_2(\Gamma(V_{8n}))}{|e(\Gamma(V_{8n}))|}-\frac{M_1(\Gamma(V_{8n}))}{|v(\Gamma(V_{8n}))|}&= \frac{144n^4+168n^3+20n^2-4n}{12n^2+4n}-\frac{80n^3+80n^2-8n}{8n} \\
        &= \frac{192n^5+64n^4-64n^3}{8n(12n^2+4n)}>0.
     \end{align*}
 
{\it Part 3.} Suppose that $\Gamma(V_{8n}) =\OSCom(V_{8n})$. If $n$ is even, then, by Theorem \ref{OSCom V_8n}, we have $\Gamma(V_{8n}) \cong K_{8n}$. As such, $M_1(\Gamma(V_{8n}))=8n(8n-1)^2$ and $M_2(\Gamma(V_{8n}))=\binom{8n}{2}(8n-1)^2=\frac{8n(8n-1)^3}{2}$. Also, $|v(\Gamma(V_{8n}))|=|V_{8n}|=8n$ and $|e(\Gamma(V_{8n}))|=\binom{8n}{2}=\frac{8n(8n-1)}{2}$. Therefore
\begin{align*}
    \frac{M_2(\Gamma(V_{8n}))}{|e(\Gamma(V_{8n}))|}=(8n-1)^2=\frac{M_1(\Gamma(V_{8n}))}{|v(\Gamma(V_{8n}))|}.
\end{align*}
If $n$ is odd, then, by Theorem \ref{OSCom V_8n}, we have $\Gamma(V_{8n})\cong K_{2n+4} \vee (K_{2n} \sqcup K_{4n-4})$. Now, by Lemma \ref{Zagreb}, we have 
\[
M_1(\Gamma(V_{8n})) = 304n^3+80n^2+8n \text{ and } M_2(\Gamma(V_{8n})) = 960n^4+1080n^3+1568n^2-220n.
\] 
      Also, $|v(\Gamma(V_{8n}))|=|V_{8n}|=8n$ and $|e(\Gamma(V_{8n}))|=\binom{2n+4}{2}+\binom{2n}{2}+\binom{4n-4}{2}+2n(2n+4)+(2n+4)(4n-4)=24n^2+4n$. Therefore
     \begin{align*}
         \frac{M_2(\Gamma(V_{8n}))}{|e(\Gamma(V_{8n}))|}-\frac{M_1(\Gamma(V_{8n}))}{|v(\Gamma(V_{8n}))|}&= \frac{240n^3+270n^2+392n-55}{6n+1)}-(38n^2+10n+1) \\
         &= \frac{12n^3+172n^2+376n-56}{6n+1}> 0.
     \end{align*}
\end{proof}
\begin{theorem}\label{Zagreb SD_8n}
The graphs $\ESCom(SD_{8n})$, $\CSCom(SD_{8n})$ and $\OSCom(SD_{8n})$ satisfy the Hansen--Vuki{\v{c}}evi{\'c} conjecture \eqref{Conj-eq}.
\end{theorem}
\begin{proof}
   {\it Part 1.} Suppose that $\Gamma(SD_{8n})=\ESCom(SD_{8n})$. As such, by Proposition \ref{ESCom(G)}(e), we have $\Gamma(SD_{8n}) \cong K_2 \vee (K_{4n-2} \sqcup 2nK_2)$ or $K_4 \vee (K_{4n-4} \sqcup nK_4)$ according as $n$ is even or odd. Hence, the result follows from Part 1 of Theorem \ref{Zagreb V_8n} for both the cases. 

    {\it Part 2.} Suppose that $\Gamma(SD_{8n})=\CSCom(SD_{8n})$. If $n$ is even,  then, by Theorem \ref{CSCom SD_8n}, we have $\Gamma(SD_{8n}) \cong K_2 \vee(K_{4n-2} \sqcup 2K_{2n})$. As such, the result follows from Part 2 of Theorem \ref{Zagreb V_8n}. If $n$ is odd,  then, by Theorem \ref{CSCom SD_8n}, we have $\Gamma(SD_{8n}) \cong K_4 \vee(K_{4n-4} \sqcup K_{4n})$. Now, by Lemma \ref{Zagreb}, we have 
     \[
     M_1(\Gamma(SD_{8n})) = 128n^3+256n^2+8n \text{ and } M_2(\Gamma(SD_{8n})) = 256n^4+832n^3+336n^2-52n.
     \]
 Also, $|v(\Gamma(SD_{8n}))|=|SD_{8n}|=8n$ and $|e(\Gamma(SD_{8n}))|=6+4(4n-4)+(2n-2)(4n-5)+2n(4n-1)+4 \times 4n=16n^2+12n$. Therefore
 \begin{align*}
     \frac{M_2(\Gamma(SD_{8n}))}{|e(\Gamma(SD_{8n}))|}-\frac{M_1(\Gamma(SD_{8n}))}{|v(\Gamma(SD_{8n}))|}= \frac{1024n^4-512n^2-512n^2}{8n(16n^2+12n)} = \frac{512n^2(2n^2-n-1)}{8n(16n^2+12n)}.
 \end{align*} 
 It is easy to see that $512n^2(2n^2-n-1)>0$ and $8n(16n^2+12n)>0$ for all $n \geq 3$ such that $n$ is odd.  

  {\it Part 3.} Suppose that $\Gamma(SD_{8n})=\OSCom(SD_{8n})$. By Theorem \ref{OSCom SD_8n}, we have $\Gamma(SD_{8n})$ $ \cong K_{8n}$ and hence the result follows from Part 3 of Theorem \ref{Zagreb V_8n}.
\end{proof}
Since $SD_{8n} = QD_{2^m}$ for $n=2^{m-3}$ we have the following corollary.
\begin{cor}
    The graphs $\ESCom(QD_{2^m})$, $\CSCom(QD_{2^m})$ and $\OSCom(QD_{2^m})$ satisfy the Hansen--Vuki{\v{c}}evi{\'c} conjecture \eqref{Conj-eq}.
\end{cor}
\begin{theorem}
   The graphs $\ESCom(U_{6n})$, $\CSCom(U_{6n})$ and $\OSCom(U_{6n})$ satisfy the Hansen--Vuki{\v{c}}evi{\'c} conjecture \eqref{Conj-eq}. 
\end{theorem}
\begin{proof}
{\it Part 1.} Suppose that $\Gamma(U_{6n})=\ESCom(U_{6n})$. By Proposition \ref{ESCom(G)}(f), we have $\Gamma(U_{6n}) \cong K_n \vee (K_{2n} \sqcup 3K_n)$. So, by Lemma \ref{Zagreb}, we have 
\[
M_1(\Gamma(U_{6n})) = 66n^3-36n^2+6n \text{ and } M_2(\Gamma(U_{6n})) = \frac{228n^4-198n^3+54n^2-6n}{2}.
\]
     Also, $|v(\Gamma(U_{6n}))|=|U_{6n}|=6n$ and $|e(\Gamma(U_{6n}))|=\binom{n}{2}+n \times 2n +\binom{2n}{2}+n \times 3n +3\binom{n}{2}=\frac{18n^2-6n}{2}$. Therefore, for all $n \geq 1$, we have
     \begin{align*}
        \frac{M_2(\Gamma(U_{6n}))}{|e(\Gamma(U_{6n}))|}-\frac{M_1(\Gamma(U_{6n}))}{|v(\Gamma(U_{6n}))|}= \frac{180n^5-144n^4}{6n(18n^2-6n)}>0.
     \end{align*}
 
     {\it Part 2.} Suppose that $\Gamma(U_{6n})=\CSCom(U_{6n})$ or $\OSCom(U_{6n})$. By Theorem \ref{CSCom U_6n} and Theorem \ref{OSCom U_6n}, we have $\Gamma(U_{6n}) \cong K_n \vee(K_{2n} \sqcup K_{3n})$. So, by Lemma \ref{Zagreb}, we have 
     \[
     M_1(\Gamma(U_{6n})) = 102n^3-48n^2+6n \text{ and } M_2(\Gamma(U_{6n})) = \frac{432n^4-306n^3+72n^2-6n}{2}.
     \]
     Also, $|v(\Gamma(U_{6n}))|=|U_{6n}|=6n$ and $|e(\Gamma(U_{6n}))|=\binom{n}{2}+n \times 2n +\binom{2n}{2}+n \times 3n +\binom{3n}{2}=\frac{19n^2-6n}{2}$. Therefore,  for all $n \geq 1$, we have
     \begin{align*}
        \frac{M_2(\Gamma(U_{6n}))}{|e(\Gamma(U_{6n}))|}-\frac{M_1(\Gamma(U_{6n}))}{|v(\Gamma(U_{6n}))|}= \frac{656n^5-312n^4+30n^3}{6n(19n^2-6n)} >0.
     \end{align*}
\end{proof}
\begin{theorem}
   The graphs $\ESCom(M_{2mn})$, $\CSCom(M_{2mn})$ and $\OSCom(M_{2mn})$ satisfy the Hansen--Vuki{\v{c}}evi{\'c} conjecture \eqref{Conj-eq}.
\end{theorem}
\begin{proof}
{\it Part 1.} Suppose that $\Gamma(M_{2mn})=\ESCom(M_{2mn})$. If $m$ is odd, then, by Proposition \ref{ESCom(G)}(g), we have $\Gamma(M_{2mn}) \cong  K_n \vee (K_{(mn-n)} \sqcup mK_n)$. Now, by Lemma \ref{Zagreb}, we have 
\[
M_1(\Gamma(M_{2mn})) = m^3n^3-2m^2n^2+2mn+3m^2n^3-6mn^2+4mn^3 \quad \text{and}
\]
\begin{align*}
 M_2(\Gamma(M_{2mn})) = \frac{1}{2}&\left(m^4n^4-3m^3n^3+9mn^2+2m^3n^4-9m^2n^3+9m^2n^4\right.\\
&\quad \quad \quad \quad \quad \quad \quad \quad \quad \quad \quad\left.-12mn^3+3m^2n^2 -2mn+4mn^4 \right).
\end{align*}
Also, $|v(\Gamma(M_{2mn}))|=|M_{2mn}|=2mn$ and $|e(\Gamma(M_{2mn}))|=\binom{n}{2}+n^2(m-1)+mn^2+\binom{mn-n}{2}+m\binom{n}{2}=\frac{m^2n^2+3mn^2-2mn}{2}$. Therefore
\begin{align*}
    \frac{M_2(\Gamma(M_{2mn}))}{|e(\Gamma(M_{2mn}))|}-\frac{M_1(\Gamma(M_{2mn}))}{|v(\Gamma(M_{2mn}))|}&=\frac{m^5n^5-2m^4n^4-2m^4n^5+5m^3n^5+2m^2n^4-4m^2n^5}{2mn(m^2n^2+3mn^2-2mn)} \\
    &= \frac{m^4n^4(mn-2n-2)+m^2n^5(5m-4)+2m^2n^4}{2mn(m^2n^2+3mn^2-2mn)}.
\end{align*}
It is easy to see that $m^4n^4(mn-2n-2)+m^2n^5(5m-4)+2m^2n^4>0$ and $2mn(m^2n^2+3mn^2-2mn)>0$ for all $n \geq 1$ and $m \geq 3$ such that $m$ is odd.

 If $m$ is even, then, by Proposition \ref{ESCom(G)}(g), we have $\Gamma(M_{2mn}) \cong K_{2n} \vee (K_{(\frac{m}{2}-1)2n} \sqcup \frac{m}{2}K_{2n})$. Now, by Lemma \ref{Zagreb}, we have 
 \[
 M_1(\Gamma(M_{2mn})) = m^3n^3-2m^2n^2+6m^2n^3+16mn^3-12mn^2+2mn \quad \text{and}
 \]
 \begin{align*}
    M_2(\Gamma(M_{2mn}))
    &=\frac{1}{2}\left(m^4n^4-3m^3n^3+18mn^2+4m^3n^4-18m^2n^3+36m^2n^4 \right. \\
    & \quad \quad \quad \quad \quad \quad \quad \quad \quad \quad \quad \quad \quad \quad \quad \left.-48mn^3+3m^2n^2-2mn+32mn^4 \right).
\end{align*}
Also, $|v(\Gamma(M_{2mn}))|=|M_{2mn}|=2mn$ and $|e(\Gamma(M_{2mn}))|=\binom{2n}{2}+2n(mn-2n)+2n \times 2n \times \frac{m}{2}+\binom{mn-2n}{2}+\frac{m}{2}\binom{2n}{2}=\frac{m^2n^2+6mn^2-2mn}{2}$. Therefore, for all $m \geq 6$ and $n \geq 1$, we have
\begin{align*}
    \frac{M_2(\Gamma(M_{2mn}))}{|e(\Gamma(M_{2mn}))|}-\frac{M_1(\Gamma(M_{2mn}))}{|v(\Gamma(M_{2mn}))|}&=\frac{m^5n^5-2m^4n^4-4m^4n^5+20m^3n^5+8m^2n^4-32m^2n^5}{2mn(m^2n^2+6mn^2-2mn)} \\
    &= \frac{m^4n^4(mn-4n-2)+m^2n^5(20m-32)+8m^2n^4}{2mn(m^2n^2+6mn^2-2mn)} >0.
\end{align*}

 {\it Part 2.} Suppose that $\Gamma(M_{2mn})=\CSCom(M_{2mn})$. If $m$ is odd, then, by Theorem \ref{CSCom M_2mn}, we have $\Gamma(M_{2mn}) \cong K_n \vee (K_{(mn-n)} \sqcup K_{mn})$. So, by Lemma \ref{Zagreb}, we have 
\[
M_1(\Gamma(M_{2mn})) = 2m^3n^3+5m^2n^3-4m^2n^2+mn^3-4mn^2+2mn \quad \text{and}
\]
\begin{align*}
     &M_2(\Gamma(M_{2mn}))\\
    &\quad=\frac{1}{2}\left(2m^4n^4-6m^3n^3+6mn^2+8m^3n^4-15m^2n^3+6m^2n^4-3mn^3+6m^2n^2 -2mn \right).
\end{align*}
Also, $|v(\Gamma(M_{2mn}))|=|M_{2mn}|=2mn$ and $|e(\Gamma(M_{2mn}))|=\binom{n}{2}+n^2(m-1)+mn^2+\binom{mn-n}{2}+\binom{mn}{2}=\frac{2m^2n^2+2mn^2-2mn}{2}$. Therefore, for all $m \geq 3$ and $n \geq 1$, we have
\begin{align*}
    \frac{M_2(\Gamma(M_{2mn}))}{|e(\Gamma(M_{2mn}))|}-\frac{M_1(\Gamma(M_{2mn}))}{|v(\Gamma(M_{2mn}))|}&=\frac{2m^4n^5-4m^3n^4+4m^2n^4-2m^2n^5}{2mn(2m^2n^2+2mn^2-2mn)} \\
    &= \frac{2m^2n^4\{(m(mn-2)-(n-2)\}}{2mn(2m^2n^2+2mn^2-2mn)} >0.
\end{align*}

 If $m$ and $\frac{m}{2}$ are even, then, by Theorem \ref{CSCom M_2mn}, we have $\Gamma(M_{2mn}) \cong K_{2n} \vee (2K_{\frac{mn}{2}} \sqcup K_{(mn-2n)})$. Now, by Lemma \ref{Zagreb}, we have 
 \[
 M_1(\Gamma(M_{2mn})) = \frac{5}{4}m^3n^3-6m^2n^2+8m^2n^3+4mn^3-8mn^2+2mn \quad \text{and}
 \]
 \begin{align*}
    M_2(\Gamma(M_{2mn}))
    =\frac{9}{16}m^4n^4-\frac{15}{8}m^3n^3+7mn^2&+\frac{9}{2}m^3n^4-12m^2n^3\\
    &+11m^2n^4-6mn^3+\frac{9}{4}m^2n^2-mn.
\end{align*}
Also, $|v(\Gamma(M_{2mn}))|=|M_{2mn}|=2mn$ and $|e(\Gamma(M_{2mn}))|=\binom{2n}{2}+2n(mn-2n)+2n \times \frac{mn}{2} \times 2 +\binom{mn-2n}{2}+2\binom{\frac{mn}{2}}{2}=\frac{3}{4}m^2n^2+2mn^2-mn$. Therefore, for all $m \geq 8$ and $n \geq 1$, we have
\begin{align*}
    &~~~~~~\frac{M_2(\Gamma(M_{2mn}))}{|e(\Gamma(M_{2mn}))|}-\frac{M_1(\Gamma(M_{2mn}))}{|v(\Gamma(M_{2mn}))|} \\
    &=\frac{\frac{3}{16}m^5n^5+4m^4n^4-3m^3n^3+2m^3n^4-\frac{1}{2}m^4n^5+3m^3n^5+8m^2n^4-8m^2n^5+2m^2n^3}{2mn(\frac{3}{4}m^2n^2+2mn^2-mn)} \\
    &= \frac{\frac{3}{16}m^5n^5+m^3n^3(4mn-3)+2m^3n^4-\frac{1}{2}m^4n^5+8m^2n^4+m^2n^5(3m-8)+2m^2n^3}{2mn(\frac{3}{4}m^2n^2+2mn^2-mn)} \\
    &>0.
\end{align*}

If $m$ is even and $\frac{m}{2}$ is odd, then, by Theorem \ref{CSCom M_2mn}, we have $\Gamma(M_{2mn}) \cong K_{2n} \vee (K_{mn} \sqcup K_{(mn-2n)})$. Now, by Lemma \ref{Zagreb}, we have 
\[
M_1(\Gamma(M_{2mn})) = 2m^3n^3+10m^2n^3-4m^2n^2+4mn^3-8mn^2+2mn \quad \text{and}
\]
\begin{align*}
     &M_2(\Gamma(M_{2mn}))\\
    & \quad=m^4n^4-3m^3n^3+6mn^2+8m^3n^4-16m^2n^3+12m^2n^4-6mn^3+3m^2n^2-mn.
\end{align*}
Also, $|v(\Gamma(M_{2mn}))|=|M_{2mn}|=2mn$ and $|e(\Gamma(M_{2mn}))|=\binom{2n}{2}+2n(mn-2n)+2n \times mn+\binom{mn-2n}{2}+\binom{mn}{2}=m^2n^2+2mn^2-mn$. Therefore, for all $m \geq 6$ and $n \geq 1$, we have
\begin{align*}
    \frac{M_2(\Gamma(M_{2mn}))}{|e(\Gamma(M_{2mn}))|}-\frac{M_1(\Gamma(M_{2mn}))}{|v(\Gamma(M_{2mn}))|}&=\frac{2m^4n^5-6m^3n^4+8m^2n^4-8m^2n^5}{2mn(m^2n^2+2mn^2-mn)} \\
    &= \frac{m^2n^4\{n(2m^2-8)-(6m-8)\}}{2mn(m^2n^2+2mn^2-mn)}>0.
\end{align*}
\end{proof}

{\bf Acknowledgement.} The first author is thankful to Council of Scientific and Industrial Research  for the fellowship (File No. 09/0796(16521)/2023-EMR-I).

\section*{Declarations}
\begin{itemize}
	\item Funding: No funding was received by the authors.
	\item Conflict of interest: The authors declare that they have no conflict of interest.
	\item Availability of data and materials: No data was used in the preparation of this manuscript.
\end{itemize}

\end{document}